\RequirePackage{fix-cm}
\documentclass{svjour3}                     
\makeatletter
\def\cl@chapter{}
\makeatother
\usepackage[T1]{fontenc}

\smartqed  
\usepackage{graphicx}
\input{definitions}

\begin{document}

\title{Generalized Parallel Tempering on Bayesian Inverse Problems\thanks{}
}


\author{Jonas Latz         \and
        Juan P. Madrigal-Cianci \and
        Fabio Nobile \and 
        Ra\'ul Tempone. 
}


\institute{J. Latz, \at
	Department of Applied Mathematics and Theoretical Physics, University of Cambridge.
	\and
	J.P Madrigal-Cianci and F. Nobile, \at
	SB-MATH-CSQI, \'{E}cole Polytechnique F\'{e}d\'{e}rale de Lausanne. \\
	\email{juan.madrigalcianci@epfl.ch}
	\and
	R. Tempone, \at
	Computer, Electrical and  Mathematical Sciences and Engineering, KAUST, and  Alexander von Humboldt professor in Mathematics of Uncertainty Quantification, RWTH Aachen University.
}
\date{Received: date / Accepted: date}

\maketitle

\begin{abstract}
In the current work we present two  generalizations of the Parallel Tempering algorithm \rjp{in the context of discrete-time Markov chain Monte Carlo methods for Bayesian inverse problems. These generalizations use state-dependent swapping rates, inspired by the so-called continuous time Infinite Swapping algorithm presented in [Plattner, Nuria, et al. 2011. J. Chem Phys.  135.13: 134111]}. We analyze the reversibility and ergodicity properties of our generalized PT algorithms. Numerical results on sampling from different target distributions, show that the proposed methods significantly improve sampling efficiency over more traditional sampling algorithms such as Random Walk Metropolis, \rjp{preconditioned Crank-Nicolson,} and  (standard)  Parallel Tempering.

This article has been published as Latz, J., Madrigal-Cianci, J.P., Nobile, F. et al. Generalized parallel tempering on Bayesian inverse problems. \textit{Stat Comput} \textbf{31}, 67 (2021). \url{https://doi.org/10.1007/s11222-021-10042-6}, and can be accessed here: \url{https://rdcu.be/cznfD}

%

\end{abstract}

\section{Introduction}\label{S:Intro}

Modern computational facilities and recent advances in computational techniques have made the use of Markov Chain Monte Carlo (MCMC) methods feasible for some large-scale Bayesian inverse problems (BIP), where the goal is to characterize the posterior distribution of a set of parameters $\te$ of a computational model which describes some physical phenomena, conditioned on some (usually indirectly) measured data $y$. However, some computational difficulties are prone to  arise when dealing with \textit{difficult to explore}  posteriors, i.e., posterior distributions that are multi-modal,  or that concentrate around  a non-linear, lower-dimensional manifold, since some of the more commonly-used Markov transition kernels in MCMC algorithms, such as random walk Metropolis (RWM) or preconditioned Crank-Nicholson (pCN), are not well-suited in such situations. This in turn can make the computational time needed to properly \textit{explore} these complicated target distributions arbitrarily long. Some recent works address these issues by employing Markov transitions kernels that use geometric information \cite{beskos2017geometric}; however, this requires  efficient computation of the gradient of the posterior density, which might not always be feasible, particularly when the underlying computational model is a so-called ``black-box''. 

In recent years, there has been an active development of computational techniques and algorithms to overcome these issues using a \emph{tempering strategy} \cite{dia2019continuation,earl2005parallel,latz2017multilevel,miasojedow2013adaptive,vrugt2009accelerating}. Of particular importance for the work presented here is the Parallel Tempering (PT) algorithm \cite{earl2005parallel,lkacki2016state,miasojedow2013adaptive} (also known as \emph{replica exchange}), which finds its origins in the physics and molecular dynamics community. The general idea behind such methods is to simultaneously run $K$ independent  MCMC chains, where each chain is invariant with respect to a \textit{flattened} (referred to as \textit{tempered}) version of the posterior of interest $\mu$, while, at the same time, proposing to swap states between any two chains every so often. Such a swap is then accepted using the standard Metropolis-Hastings (MH) acceptance-rejection rule. Intuitively, chains with a larger smoothing parameter (referred to as \textit{temperature}) will  be able to better explore the parameter space. Thus, by proposing to exchange states between chains that target posteriors at different temperatures, it is possible for the chain of interest (i.e., the one targeting $\mu$) to mix faster, and to avoid the undesirable behavior of some MCMC samplers of getting  ``stuck'' in a mode. Moreover, the fact that such an exchange of states is accepted with the typical MH acceptance-rejection rule, will guarantee that the chain targeting $\mu$ remains invariant with respect to such probability measure \cite{earl2005parallel}. 

Tempering ideas have been successfully used to sample from posterior distributions arising in different fields of science, ranging from astrophysics to machine learning \cite{desjardins2010tempered,earl2005parallel,miasojedow2013adaptive,van2008parameter}.  The works \cite{madras2002markov,woodard2009conditions} have studied the convergence of the PT algorithm  from a theoretical perspective and provided minimal conditions for its rapid mixing. 
Moreover, the idea of tempered distributions has not only been applied in combination with parallel chains. For example,
the simulated tempering method \cite{marinari_1992} uses a single chain and varies the temperature within this chain. In addition, tempering forms the basis of efficient particle filtering methods for stationary model parameters in Sequential Monte Carlo settings \cite{Beskos16,Beskos15,kahle2018bayesian,Kantas14,latz2017multilevel} and Ensemble Kalman Inversion \cite{Schillings17}.

A generalization over the PT approach, originating from the molecular dynamics community, is the so-called \emph{Infinite Swapping (IS)} algorithm \cite{dupuis2012infinite,plattner2011infinite}. 
As opposed to PT, this IS paradigm is a continuous-time Markov process and  considers the limit where states between chains are swapped infinitely often. 
It is shown in \cite{dupuis2012infinite} that such an approach can in turn be understood as a swap of dynamics, i.e., kernel and temperature (as opposed to states) between chains. We remark that once such a change in dynamics is considered, it is not possible to distinguish particles belonging to different chains. However, since the stationary distribution of each chain is known, importance sampling can be employed to compute posterior expectations with respect to the target measure of interest.
Infinite Swapping has been successfully applied in the context of computational molecular dynamics and rare event simulation \cite{doll2012rare,lu2013infinite,plattner2011infinite}; however, to the best of our knowledge, such methods have not been implemented in the context of Bayesian inverse problems. 

In light of this, the current work aims at importing such ideas to the BIP setting, 
by presenting them in a discrete-time Metropolis-Hastings Markov chain Monte Carlo context. We will refer to these algorithms as \emph{Generalized Parallel Tempering} (GPT).  We emphasize, however, that these methods are \emph{not} a time discretization of the continuous-time Infinite Swapping presented in \cite{dupuis2012infinite}, but, in fact, a discrete-time Markov process inspired by the ideas presented therein with suitably defined state-dependent probabilities of swapping states or dynamics. We now summarize the main contributions of this work.
	
\rjp{We begin by presenting a generalized framework for discrete time PT in the context of MCMC for BIP, and then proceed to propose, analyze and implement two novel state-dependent PT algorithms inspired by the ideas presented in \cite{dupuis2012infinite}}.

 Furthermore, we prove that our GPT methods  have the right invariant measure, by showing reversibility of the generated Markov chains, and prove their ergodicity. Finally, we implement the proposed GPT algorithms for an array of  Bayesian inverse problems, \rjp{comparing their efficiency to that of an un-tempered, (single temperature), version of the underlying MCMC algorithm, and standard PT. For the base method to sample at the cold temperature level, we use Random Walk Metropolis (RWM) (Sections \ref{ss:exp_circle}-\ref{ss:exp_bip_ac}) or preconditioned Crank-Nicolson (Section \ref{ss:high_dim}), however, we emphasize that our methods can be used together with any other, more efficient base sampler.} Experimental results show improvements in terms of computational efficiency of GPT over un-tempered RWM and standard PT,  thus making the proposed methods  attractive from  a computational perspective. \rjp{From an implementation perspective, the swapping component of  our proposed methods is rejection-free, thus effectively eliminating some tuning parameters on the PT algorithm, such as swapping frequency.} 
 
 We notice that a PT algorithm with state-dependent swapping probabilities has been proposed in \cite{lkacki2016state}, however, such a work only consider pairwise swapping of chains and a different construction of the swapping probabilities, resulting in a less-efficient sampler, at least for the BIPs addressed in this work.

 Our ergodicity result relies on an $L_2$ spectral gap analysis.  It is known  \cite{rudolf2011explicit} that when a Markov chain is both reversible and has a positive $L_2$-spectral gap, one can in turn provide non-asymptotic error bounds on the mean square error of an ergodic estimator of the chain. Our bounds on the $L_2$-spectral gap, however, are far from being sharp and could possibly be improved using e.g., domain decomposition ideas as in \cite{woodard2009conditions}. Such analysis is left for a future work.

The rest of this paper is organized as follows. Section \ref{S:Setting} is devoted to the introduction of the notation, Bayesian inverse problems, and Markov chain Monte Carlo methods. In Section \ref{S:pt_e_is}, we provide a brief review of (standard) PT (Section \ref{ss:pt}), and introduce the two versions of the GPT algorithm in Sections \ref{SS:uwis} and \ref{SS:wis}, respectively. In fact, we present a general framework that accommodates both the standard PT algorithms and our generalized versions. In Section \ref{S:AsympInfSwap}, we recall some of the standard theory of Markov chains in Section \ref{ss:prelim} and state the main theoretical results of the current work (Theorems \ref{proposition:uwIS_ergodicity} and \ref{proposition:wIS_ergodicity}) in Section \ref{ss:main_results}. The proof of these Theorems is given by a series of Propositions in Section \ref{ss:main_results}. We present some numerical experiments in Section \ref{S:Num_Exp}, and draw some conclusions in Section \ref{S:conclusions}. 

\color{black}

\section{Problem setting}\label{S:Setting}
\subsection{Notation}\label{S:BIP}

 Let $(W,\lno \cdot \rno)$ be a separable Banach space with associated Borel $\sigma$-algebra $\mathcal{B}(W)$, and let $\nu_W$ be a 
 $\sigma$-finite ``reference'' measure on $W$.   For any measure $\mu$ on $(W, \mathcal{B}(W))$ that is absolutely continuous with respect to $\nu_W$ (in short $\mu \ll \nu_W$), we define the Radon-Nikodym derivative $\pi_\mu := \frac{\mathrm{d}\mu}{\mathrm{d}\nu_W}$. We denote by $\overline{\mathcal{M}}(W)$ the set of real-valued signed measures on $(W,\BW)$, and by $\mathcal{M}(W)\subset \overline{\mathcal{M}}(W)$ the set of probability measures on $(W,\BW)$. 	
%
%

Let  $W_1, W_2$ be two separable Banach spaces with reference measures $\nu_{W_1},\nu_{W_2}$, and let $\mu_1 \ll \nu_{W_1}, \mu_{2} \ll \nu_{W_2}$ be two probability measures, with corresponding densities  given by $\pi_{\mu_1},\pi_{\mu_2}$.
The \emph{product} of these two measures is defined by
\begin{align}
\bm \mu(A)&=\left(\mu_1 \times \mu_{2}\right)(A)\\ & = \iint_A \pi_{\mu_1}(\te_1) \pi_{\mu_2}(\te_2)\nu_{W_1}(\mathrm{d}\te_1)\nu_{W_2}(\mathrm{d}\te_2), 
\end{align}
for all $A\in \mathcal{B}(W_1 \times W_2).$ Joint measures on $(W_1\times W_2, \mathcal{B}(W_1,\times W_2))$ will always be written in boldface. 

A \emph{Markov kernel} on a Banach space $W$ is a function $p:W\times\mathcal{B}(W)\rightarrow[0,1]$ such that \begin{enumerate}
		\item For each $A$ in $\mathcal{B}(W)$, the mapping $ W\ni\te\mapsto p(\te,A)$, is a $\mathcal{B}(W)$-measurable real-valued function. 
		\item For each $\te$ in $W$, the mapping $\mathcal{B}(W)\ni A\mapsto p(\te,A)$, is a probability measure on $(W,\mathcal{B}(W))$. 	
	\end{enumerate}

 We denote by $P$ the Markov transition operator associated to $p$, \rjp{which acts on measures as $\nu\mapsto\nu P\in\mathcal{M}(W),$ and on functions as $f\mapsto Pf,$ with $Pf$ measurable on $(W,\mathcal{B}(W)),$} such that 
	\begin{align}
(\nu P) (A)&=\int_W p(\te,A)\nu(\mathrm{d}\te), \quad  \forall A\in  \BW,\\
(Pf)(\te)&=\int_W f(z)p(\te,\mathrm{d}z), \quad \forall \te \in W .
	\end{align}

Let  $P_k,\ k=1,2,$ be Markov  transition operators  associated to kernels $p_k:W_k\times\mathcal{B}(W_k)\mapsto [0,1].$ We define the \emph{tensor product Markov operator} $\mathbf{P}:=P_1\otimes P_2$ as the Markov operator associated with the product measure $\mathbf{p}(\teb,\cdot)=p_1(\te_1,\cdot)\times p_2(\te_2,\cdot), \ \teb=(\te_1,\te_2)\in W_1\times W_2$. In particular, $\bm \nu  \mathbf{P}$ is the measure on $(W_1 \times W_2, \mathcal{B}(W_1 \times W_2))$ that satisfies
	 \begin{align}
(\bm \nu  \mathbf{P}) (A_1 \times A_2)&=\iint_{W_1 \times W_2} p_1(\te_1,A_1)p_2(\te_2,A_2)\bm\nu(\mathrm{d}\theta_1, \mathrm{d}\theta_2),
\end{align}
for all $A_1 \in \mathcal{B}(W_1)$ and  $A_2 \in \mathcal{B}(W_2)$.
	Moreover, $(\mathbf{P}f): W_1 \times W_2 \rightarrow \mathbb{R}$ is the function given by
	\begin{align}
(\mathbf{P}f)(\teb)&=\iint_{W_1 \times W_2}f(z_1,z_2)p_1(\te_1,\mathrm{d}z_1)p_2(\te_2,\mathrm{d}z_2) ,
	\end{align}
	for an appropriate $f: W_1 \times W_2 \rightarrow \mathbb{R}$.
	
%


We say that a Markov operator $P$ (resp. $\mathbf{P}$) is \emph{invariant} with respect to a measure $\nu$ (resp. $\bm \nu$) if $\nu P=\nu$ (resp.  $\bm\nu \mathbf{P}=\bm\nu$ ). A related concept to invariance is that of reversibility;	a Markov kernel $p:W\times \BW\mapsto[0,1]$ is said to be \emph{reversible}  (or \emph{$\nu$-reversible}) with respect to  a measure  $\nu\in \mathcal{M}(W)$ if  \begin{align}\label{Eq:Rev_1}\int_B p(\te,A)\nu(\mathrm{d}\te)=\int_A p(\te,B)\nu(\mathrm{d}\te),  \quad  \forall A,B\in \BW.\end{align}

 For two given $\nu$-invariant Markov operators $P_1,P_2$, we say that $P_1P_2$ is a \emph{composition} of Markov operators, not to be confused with $ P_1\otimes P_2$. Furthermore, given a composition of $K$ $\nu$-invariant Markov operators $P_c:=P_1P_2\dots P_K$, we say that $P_c$ is \emph{palindromic} if $P_1=P_K,$ $P_2=P_{K-1},$ \dots, $P_k=P_{K-{k+1}}, \ k=1,2\dots,K$. It is known (see, e.g., \cite[Section 1.12.17]{brooks2011handbook}) that a palindromic, $\nu$-invariant Markov operator $P_c$ has an associated Markov transition  kernel $p_c$ which is $\nu$-reversible.

\subsection{Bayesian inverse problems}\label{Ss:Bckgnd}

 Let $(\Theta,\lno\cdot\rno_\Theta)$ and $(Y,\lno\cdot\rno_Y)$ be separable Banach spaces with associated Borel $\sigma$-algebras $\mathcal{B}(\Theta),\ \mathcal{B}(Y)$. In Bayesian Inverse Problems we aim at characterizing the probability distribution of a set of parameters $\te\in\Theta$ conditioned on some measured data $y\in Y$, where the relation between $
	\te$ and $y$ is given by:
\begin{align} \label{Eq:InvProb}
{y}={\eff}(\te)+{\varepsilon}, \quad \varepsilon\sim \mu_\textrm{noise}.
\end{align}
Here $\varepsilon$ is some random noise with known distribution $\mu_\mathrm{noise}$ (assumed to have a density $\pi_\textrm{noise}$ with respect to some reference measure $\nu_Y$ on $Y$) and $\mathcal{F}:\Theta\mapsto Y$ is the so-called \emph{forward mapping operator}. Such an operator can model, e.g., the numerical solution of a possibly non-linear Partial Differential Equation (PDE) which takes $
\te$ as a set of parameters. \rjp{We will assume that the data $y$ is finite dimensional, i.e., $Y=\R^{d_y}$ for some $d_y\geq1$}, and that $\te\sim \mu_\mathrm{pr}$. Furthermore, we define the \emph{potential} function $\Phi(\te;y):\Theta\times Y\mapsto \R$ as
\begin{align}\label{eq:potential}
\Phi(\te;y)=-\log \left[\pi_\textrm{noise}(y-\mathcal{F}(\te))\right],
\end{align}
where the function $\Phi(\te;y)$ is a measure of the misfit between the recorded data $y$ and the predicted value $\mathcal{F}(\te)$, and often depends only on $\lno y-\mathcal{F}(\te)\rno_Y$. The
 solution to the Bayesian inverse problem is given by  (see, for example, \cite[Theorem 2.5]{Latz19Onthe}) \begin{align}\label{Eq:BT} {\pi(\te)}:=\pi(\te|\yy)=\frac{1}{Z}{e^{-\pot}}{\pi_\mathrm{pr}(\te)}, 
\end{align} 
where $\mu$ (with corresponding $\nu_\Theta$-density $\pi$) is referred to as the \textit{posterior measure} and $$Z:= \int_\Theta \exp(-\pot)\mu_\textrm{pr}(\mathrm{d}\te).$$ The Bayesian approach to the inverse problem consists of updating our knowledge concerning the parameter $\theta$, i.e., the prior, given the information that we observed in Equation \eqref{Eq:InvProb}. One way of doing so is to generate samples from the posterior measure $\mu^y$.  A common method for performing such a task is to use, for example, the Metropolis-Hastings algorithm, as detailed in the next section.  Once samples $\{\te^{(n)}\}_{n=0}^N$  have been obtained, the posterior expectation $\E_{\mu^y}[\QoI]$ of some $\mu^y$-integrable quantity of interest $\QoI:\Theta\mapsto\R$ can be approximated by the following ergodic estimator 
\begin{align}\label{eq:mc_estimator}
\E_{\mu^y}[\mathcal{Q}]\approx \widehat{\mathcal{Q}}:=\frac{1}{N-b} \sum_{n=b}^N \mathcal{Q}(\te^{(n)}),
\end{align}
where $b<N$ is the so-called \textit{burn-in} period, used to reduce the bias typically associated to MCMC algorithms.

\subsection{Metropolis-Hastings and tempering}\label{ss:MCMC_and_pt}

We briefly recall the Metropolis-Hastings algorithm \cite{hastings1970monte,metropolis1953equation}. Let $q_\mathrm{prop}:\Theta\times\mathcal{B}(\Theta)\mapsto[0,1]$ be an auxiliary kernel.  For $n=1,2,\dots,$   a candidate state $\te^*$ is sampled from $q_\mathrm{prop}(\te^n,\cdot)$, and proposed as the new state of the chain at step ${n+1}$. Such a state is then  accepted (i.e., we set  $\te^{n+1}=\te^*)$, with probability $\alpha_\textrm{MH}$, 
\begin{align}
\alpha_\text{MH}(\te^n,\te^*)=\min \left\{ 1, \frac{\pi_\mu(\te^*)q_\mathrm{prop}(\te^*,\te^n)}{\pi_\mu(\te^n)q_\mathrm{prop}(\te^n,\te^*)}  \right\},
\end{align}otherwise the current state is retained, i.e., $\te^{n+1}=\te^n$. Notice that, with  a slight abuse of notation, we are denoting kernel and density by the same symbol $q_\mathrm{prop}$. 
%
%
The Metropolis-Hastings algorithm induces the \textit{Markov transition kernel} $p:\Theta\times\mathcal{B}(\Theta)\mapsto[0,1]$ 
\begin{align}
p(\te,A)&= \int_A\alpha_\mathrm{MH}(\te,\te^*)q_\mathrm{prop}(\te,\mathrm{d\te^*}) \\ &+ \delta_{\te}(A)\int_\Theta(1-\alpha_\mathrm{MH}(\te,\te^*))q_\mathrm{prop}(\te,\mathrm{d\te^*}),
\end{align}
for every  $\te\in \Theta$ and $A\in \mathcal{B}(\Theta)$.
In most practical algorithms, the proposal state $\te^*$ is sampled from a state-dependent auxiliary kernel $q_\mathrm{prop}(\te^n,\cdot)$. Such is the case for \emph{random walk Metropolis} or \emph{preconditioned Crank Nicolson}, where $q_\mathrm{prop}(\te^n,\cdot)=\mathcal{N}(\te^n,\Sigma)$ or $q_\mathrm{prop}(\te^n,\cdot)=\mathcal{N}(\sqrt{1-\rho^2}\te^n,\rho\Sigma),$ $0<\rho<1$, respectively. However, these types of \textit{localized} proposals tend to present some undesirable behaviors when sampling from certain \textit{difficult} measures, which are, for example, concentrated over a manifold or are multi-modal \cite{earl2005parallel}. In the first case, in order to avoid a large rejection rate, the ``step-size'' $\lno \Sigma^{1/2}\rno$ of the proposal kernel must be quite small, which will in turn produce highly-correlated samples. In the second case, chains generated by these \textit{localized} kernels tend to get stuck in one of the modes. In either of these cases, very long chains are required to properly explore the parameter space.

 One way of overcoming such difficulties is to introduce tempering.  Let $\mu_k, \mu_\textrm{pr}$ be probability measures on $(\Theta,\mathcal{B}(\Theta)),$ $k=1,\dots,K,$ such that all $\mu_k$ are absolutely continuous with respect to $\mu_\textrm{pr}$, and let $\{T_k\}_{k=1}^K$ be a set of $K$ \textit{temperatures} such that $1=T_1<T_2<\dots<T_K\leq \infty$. In a Bayesian setting, $\mu_\textrm{pr}$ corresponds to the  prior measure and $\mu_k, k=1,\dots,K$ correspond to posterior measures associated to different temperatures. Denoting by $\pi_k$ the $\mu_\mathrm{pr}$-density of $\mu_k$, we set 
\begin{align}
\pi_k(\theta)&:=\frac{e^{-\Phi(\te;y)/T_k}}{Z_k},\quad \te\in \Theta,\label{Eq:RN_single}
\end{align}
 where $\quad Z_k:=\int_\Theta e^{-\Phi(\te;y)/T_k} \mu_\textrm{pr}(\mathrm{d}\te) ,$ and with $\Phi(\te;y)$ as the potential function defined in (\ref{eq:potential}). In the case where  $T_K=\infty$, we set $\mu_K=\mu_\textrm{pr}$. Notice that $\mu_1$ corresponds to the target posterior measure.

 We say that for $k=2,\dots, K,$ each measure $\mu_k$ is a  \textit{tempered} version of $\mu_1$. In general, the ${1/T_k}$ term in (\ref{Eq:RN_single}) serves as a ``smoothing''\footnote{Here, smoothing is to be understood in the sense that it \emph{flattens} the density.} factor, which in turn makes $\mu_k$ easier to explore as $T_k\to\infty$. In PT MCMC algorithms, we sample from all posterior measures $\mu_k$ simultaneously. Here, we first use a $\mu_k$-reversible Markov transition kernel  $p_k$ on each chain, and then, we propose to exchange states between chains at two consecutive temperatures, i.e., chains targeting $\mu_k,\mu_{k+1}, \ k\in\{1,\dots,K-1\}$. Such a proposed swap is then accepted or rejected with a standard Metropolis-Hastings acceptance rejection step. This procedure is presented in  Algorithm \ref{alg:ptg}. We remark that such an algorithm can be modified to, for example, propose to swap states every $N_s$ steps of the chain, or to swaps states between two chains $\mu_i,\mu_j$, with $i,j$ chosen randomly and uniformly from the index set $\{1,2,\dots,K\}$. In the next section we present the generalized PT algorithms which swap states according to a random permutation of the indices drawn from a state dependent probability.

\begin{algorithm}[tbh]
	\caption{Standard PT.}\label{alg:ptg}
	\begin{algorithmic}
		\Function{Standard PT}{$N,\{p_k\}_{k=1}^K$, $\{\pi_k\}_{k=1}^K,\  \mu_\textrm{pr}$}
		\State Sample $\te_k^{(1)}\sim \mu_\textrm{pr}, \ k=1,\dots,K$
		\State \color{gray} \# Do one step of MH on each chain \color{black}
		\For {$n=1,2,\dots,N-1$}
		\For {$k=1,\dots,K$}
		\State Sample $\te_k^{(n+1)}\sim p_k(\te_k^{(n)},\cdot)$
		\EndFor
		\State \color{gray} \# Swap states \color{black}
		\For {$k=1,2,\dots,K-1$}
		\State Swap states $\te^{(n+1)}_k$ and $\te^{(n+1)}_{k+1}$  with probability
		$$ \alpha_\mathrm{swap}=\min \left\{ 1, \frac{\pi_{k}(\te^{(n+1)}_{k+1} )\pi_{k+1}(\te^{(n+1)}_{k} )}{\pi_k(\te^{(n+1)}_{k} )\pi_{k+1}(\te^{(n+1)}_{k+1} )} \right\} $$
	
		\color{black}
		\EndFor
		\EndFor
		\State Output $\{\theta^{(n)}_1\}_{n=1}^{N}$
		\EndFunction
	\end{algorithmic}
	
\end{algorithm}


\section{Generalizing Parallel Tempering}\label{S:pt_e_is}
 Infinite Swapping was initially developed in the context of continuous-time MCMC algorithms, which were used for molecular dynamics simulations. In continuous-time PT, the swapping of the states is controlled by a Poisson process on the set $\{1,\ldots,K\}$. Infinite Swapping is the limiting algorithm obtained by letting the waiting times of this Poisson process go to zero. 
Hence, we swap the states of the chain infinitely often over a finite time interval.
We refer to  \cite{dupuis2012infinite}  for a thorough introduction and review of Infinite Swapping in continuous-time.
In Section 5 of the same article, the idea to use Infinite Swapping in time-discrete Markov chains was briefly discussed.
Inspired by this discussion, we present two  Generalizations of the (discrete-time) Parallel Tempering strategies. To that end, we propose to either (i) swap states in the chains at every iteration of the algorithm in such a way that the swap is accepted with probability one, which we will refer to as the \textit{Unweighted Generalized Parallel Tempering (UGPT)}, or (ii), swap dynamics (i.e., swap kernels and temperatures between chains) at every step of the algorithm. In this case, importance sampling must also be used when computing posterior expectations since this in turn provides a Markov chain whose invariant measure is not the desired one. We refer to this approach as \textit{Weighted Generalized Parallel Tempering (WGPT)}. We begin by introducing a common framework to both PT and the two versions of GPT. 

Let $(\Theta,\lno \cdot \rno_{\Theta})$ be a separable Banach space with associated Borel $\sigma$-algebra $\mathcal{B}(\Theta)$. Let us define the $K$-fold product space $\Theta^K := \bigtimes^K_{k=1} \Theta$,
with associated product $\sigma$-algebra $\mathcal{B}^K:= \bigotimes^K_{k=1} \mathcal{B}(\Theta)$,
as well as the product measures on $(\Theta^K,\mathcal{B}^K)$ \begin{align}
\bm \mu:=\bigtimes^K_{k=1}\mu_k, \label{Eq:Product_measure}
\end{align} where $\mu_k$ $k=1,\dots, K$ are the tempered measures with temperatures $1= T_1<T_2<T_3<\dots< T_K\leq \infty$ introduced in the previous section. Similarly, we define the product prior measure $\mub_\mathrm{pr}:=\bigtimes_{k=1}^K\mu_\mathrm{pr}.$ Notice that $\bm \mu$ has a density $\bm \pi(\teb)$ with respect to $\bm \mu_\textrm{prior}$ given by 
\begin{align}
&\bm \pi(\bm \te)=\prod_{k=1}^{K}\pi_k(\te_k)\label{Eq:RN_joint}, \quad \teb:=(\te_1,\dots,\te_K)\in \Theta^K,
\end{align}
with $\pi_i(\te)$ added subscript given as in \eqref{Eq:RN_single}.
The idea behind the tempering methods presented herein is to sample from $\bm \mu$ (as opposed to solely sampling from $\mu_1$) by creating a Markov chain obtained from the successive application of two $\bm \mu$-invariant Markov kernels $\mathbf{p}$ and $\mathbf{q}$, to some initial distribution $\bm \nu$, usually chosen to be the prior $\mub_\mathrm{pr}$. Each kernel acts as follows. Given the current state  $\teb^n=(\te^n_1,\dots,\te^n_K),$ the kernel $\mathbf{p}$, which we will call the \textit{standard MCMC kernel}, proposes a new, intermediate state $\tilde\teb^{n+1}=(\tilde\te^{n+1}_1,\dots,\tilde\te^{n+1}_K),$ possibly following the Metropolis-Hastings algorithm (or any other algorithm that generates a $\mub$-invariant Markov operator). The Markov kernel $\mathbf{p}$ is a product kernel, meaning that each component $\tilde \te_k^n$, $k=1\dots,K,$ is generated independently of the others. Then, the \textit{swapping kernel} $\mathbf{q}$ proposes a new state $\teb^{n+1}=(\te_1^{n+1},\dots,\te_K^{n+1})$ by introducing an 	``interaction'' between the components  of $\tilde \teb^{(n+1)}.$  This interaction step can be achieved, e.g., in the case of PT, by proposing to swap two components at two consecutive temperatures, i.e., components $k$ and $k+1$, and accepting this swap with a certain probability given by the usual Metropolis-Hastings acceptance-rejection rule. In general, the  swapping kernel is  applied every $N_s$ steps of the chain, for some $N_s\geq1$. 
 We will devote the following subsection to the construction of the swapping kernel  $\mathbf{q}$.  

\subsection{The swapping kernel $\mathbf{q}$}\label{ss:swapping_kernel} 
 Define $\mathscr{S}_K$ as the collection of all the bijective maps from $\{1,2,\dots,K \}$ to itself, i.e., the set of all $K!$ possible permutations of $\mathrm{id}:=\{1,\dots,K\}$.  Let $\sigma\in \mathscr{S}_K$ be a permutation,  and define the swapped state $\teb_\sigma:=(\te_{\sigma(1)},\dots,\te_{\sigma(K)}),$ and the inverse permutation $\sigma^{-1}\in \mathscr{S}_K$ such that $\sigma\circ\sigma^{-1}=\sigma^{-1}\circ\sigma=\mathrm{id}$. In addition, let $S_K\subseteq\mathscr{S}_K$ be any subset of $\mathscr{S}_K$ closed with respect to inversion, i.e., $\sigma \in S_K\implies \sigma^{-1}\in S_K$. We denote the cardinality of $S_K$ by $|S_K|.$
 
 \rjp{\begin{example}As a simple example, consider a Standard PT as in Algorithm \ref{alg:ptg} with $K=4$. In this case, we attempt to swap two contiguous temperatures $T_i$ and $T_{i+1}, \ i=1,2,3$. Thus, $S_K$ is the set of  permutations $\{\sigma_{1,2},\sigma_{2,3},\sigma_{3,4}\}$ with:
\begin{align}
\sigma_{1,2}&=(2,1,3,4), 	\\
\sigma_{2,3}&=(1,3,2,4), 	\\
\sigma_{3,4}&=(1,2,4,3).
\end{align}
Notice that each permutation is its own inverse; for example: $$\sigma_{1,2}(\sigma_{1,2})=\sigma_{1,2}((2,1,3,4))=(1,2,3,4)=\mathrm{id}.$$
\end{example}
}
To define the swapping kernel $\mathbf{q}$, we first need to define the swapping ratio and swapping acceptance probability.

 \begin{definition}[Swapping ratio]\label{def:s_r} We say that a function $r:\Theta^K\times S_K \mapsto [0,1]$ is a \emph{swapping ratio} if it satisfies the following two conditions: \begin{enumerate}
 		\item $\forall \teb\in \Theta^K$, $r(\teb,\cdot)$ is a probability mass function on $S_K$.
 		\item $\forall \sigma \in S_K$, $r(\cdot,\sigma)$ is measurable on $(\Theta^K,\mathcal{B}^K)$.
 	\end{enumerate}
 \end{definition}

\begin{definition}[Swapping acceptance probability]\label{def:s_a_p}

Let $\teb\in \Theta^K$ and $\sigma,\ \sigma^{-1}\in S_K$. We call \emph{swapping acceptance probability} the function $\alphaqt:\Theta^K\times S_K\mapsto[0,1] $ defined as \begin{equation}\label{Eq:swapping_acceptance_probability}
\alphaqt(\teb,\sigma)=\begin{cases}
\min\left\{ 1,\frac{\pib(\teb_\sigma )r(\teb_\sigma,\sigma^{-1})}{\pib(\teb)r(\teb,\sigma)} \right\}, & \text{if }r(\teb,\sigma)>0,\\
0 & \text{if }r(\teb,\sigma)=0.
\end{cases}
\end{equation}

\end{definition}
We can now define the swapping kernel $\mathbf{q}$. 
\begin{definition}[Swapping kernel]\label{def:s_k}	
Given a swapping ratio $r:\Theta^K\times S_K\mapsto[0,1]$ and its associated swapping acceptance probability  $\alpha_\mathrm{swap}:\Theta^K\times S_K\mapsto[0,1]$, we define the \emph{swapping Markov kernel} $\mathbf{q}:\Theta^K\times \mathcal{B}^K\mapsto[0,1]$ as
 \begin{align}\label{Eq:swapping_kernel}
\mathbf{q} (\teb, B)=\sum_{\sigma \in S_K}r(\teb,\sigma)\left[(1-\alphaqt(\teb,\sigma))\delta_{\teb} (B)\right. \\ \left.+\alphaqt(\teb,\sigma)\delta_{\teb_\sigma}(B)  \right], \qquad \teb\in \Theta^K,\  B\in \mathcal{B}^K, 
	\end{align}
	where $\delta_{\teb}(B)$ denotes the Dirac measure in $\teb$, i.e., $\delta_{\teb}(B)=1$ if $\teb\in B$ and 0 otherwise.
\end{definition}

The swapping mechanism should be understood in the following way: given a current state of the chain $\teb\in \Theta^K,$ the swapping kernel samples a permutation $\sigma$ from $S_K$ with probability $r(\teb,\sigma)$ and generates $\teb_\sigma.$ This permuted state is then accepted as the new state of the chain with probability $\alphaq(\teb,\sigma)$.  Notice that the swapping kernel follows a Metropolis-Hastings-like procedure with ``proposal'' distribution $r(\teb,\sigma)$  and acceptance probability $\alphaq(\teb,\sigma)$. Moreover, as detailed in the next proposition, such a kernel is reversible with respect to $\mub$, since it is a Metropolis-Hastings type kernel. 


%
%
%
\begin{proposition} \label{Prop_MK_Swap}
	The Markov kernel $\mathbf{q} $	defined in (\ref{Eq:swapping_kernel}) is reversible with respect to the product measure $\mub$ defined in (\ref{Eq:Product_measure}).

		\begin{proof}

	Let $A,B\in \mathcal{B}^K$. We want to show that $$\int_A q(\teb,B)\mub(\mathrm{d}\teb)=\int_B q(\teb,A)\mub(\mathrm{d}\teb).$$
	Thus,
	\begin{align}
	&\int_A q(\teb,B)\mub(\mathrm{d}\teb)=
	\underbrace{\sum_{\sigma\in S_K}\int_A r(\teb,\sigma)\alpha_\textrm{swap}(\teb,\sigma)\delta_{\teb_\sigma}(B)\pib(\teb)\mub_\mathrm{pr}(\mathrm{d}\teb)}_{I} \\&+
	\underbrace{\sum_{\sigma\in S_K}\int_A r(\teb,\sigma)\left(1-\alpha_\textrm{swap}(\teb,\sigma)\right)\delta_{\teb}(B)\pib(\teb)\mub_\mathrm{pr}(\mathrm{d}\teb) }_{II}.
	\end{align}
	Let $A_\sigma:=\{\bm z \in \Theta^K : \bm z_{\sigma^{-1}}\in A\}$, and, for notational simplicity, write $\min\{a,b\}=\{a\wedge b\}, \ a,b\in\R$. From $I$, we have: 
	\begin{align}
	&I=\sum_{\sigma\in S_K}\int_A  \left\{1\wedge\frac{\pib(\teb_\sigma)r(\teb_\sigma,\sigma^{-1})}{\pib(\teb)r(\teb,\sigma)}\right\} r(\teb,\sigma)\pib(\teb)\delta_{\teb_\sigma}(B)\mub_\mathrm{pr}(\mathrm{d}\teb)\\
	&=\sum_{\sigma\in S_K}\int_A  \left\{1\wedge\frac{\pib(\teb)r(\teb,\sigma)}{\pib(\teb_\sigma)r(\teb_{\sigma},\sigma^{-1})}\right\}r(\teb_\sigma,\sigma^{-1})\pib(\teb_\sigma)\delta_{\teb_\sigma}(B)\mub_\mathrm{pr}(\mathrm{d}\teb).
	\end{align}
\rjp{Then, noticing that $\mub_\mathrm{pr}$ is permutation invariant, we get} \begin{align} 
	&I=\sum_{\sigma\in S_K}\int_{A_\sigma}
	 \left\{1\wedge\frac{\pib(\teb_{\sigma^{-1}})r(\teb_{\sigma^{-1}},\sigma)}{\pib(\teb)r(\teb,\sigma^{-1})}\right\}\\ &  \qquad \qquad \times r(\teb,\sigma^{-1})\pib(\teb)\delta_{\teb}(B)\mub_\mathrm{pr}(\mathrm{d}\teb)\\
	&=\sum_{\sigma\in S_K}\int_{A_\sigma\cap B}
\left\{1\wedge\frac{\pib(\teb_{\sigma^{-1}})r(\teb_{\sigma^{-1}},\sigma)}{\pib(\teb)r(\teb,\sigma^{-1})}\right\}\\ &  \qquad \qquad \times r(\teb,\sigma^{-1})\pib(\teb)\delta_{\teb}(B)\mub_\mathrm{pr}(\mathrm{d}\teb)
\\
	&=\sum_{\sigma\in S_K}\int_{B}
\left\{1\wedge\frac{\pib(\teb_{\sigma^{-1}})r(\teb_{\sigma^{-1}},\sigma)}{\pib(\teb)r(\teb,\sigma^{-1})}\right\}\\ &  \qquad \qquad \times r(\teb,\sigma^{-1})\pib(\teb)\delta_{\teb}(A_\sigma)\mub_\mathrm{pr}(\mathrm{d}\teb)\\
	&=\sum_{\sigma\in S_K}\int_{B}
\left\{1\wedge\frac{\pib(\teb_{\sigma^{-1}})r(\teb_{\sigma^{-1}},\sigma)}{\pib(\teb)r(\teb,\sigma^{-1})}\right\}\\ &  \qquad \qquad \times r(\teb,\sigma^{-1})\pib(\teb)\delta_{\teb_{\sigma^{-1}}}(A)\mub_\mathrm{pr}(\mathrm{d}\teb)
	\end{align}
\begin{align}
&=\sum_{\sigma\in S_K}\int_{B} r(\teb,\sigma^{-1})\pib(\teb)\alpha_\textrm{swap}(\teb,\sigma^{-1})\delta_{\teb_{\sigma^{-1}}}(A)\mub_\mathrm{pr}(\mathrm{d}\teb)\\
	&=\sum_{\sigma\in S_K}\int_B r(\teb,\sigma)\pib(\teb)\alpha_\mathrm{swap}(\teb,\sigma)\delta_{\teb_\sigma}(A)\mub_\mathrm{pr}(\mathrm{d}\teb).
	\end{align}
	For the second term $II$ we simply have 
	\begin{align}
	&II=\sum_{\sigma\in S_K}\int_A r(\teb,\sigma)(1-\alpha_\mathrm{swap}(\teb,\sigma))\delta_{\teb}(B)\pib(\teb)\mub_\mathrm{pr}(\mathrm{d}\teb)\\
	&=\sum_{\sigma\in S_K}\int_{A\cap B} r(\teb,\sigma)(1-\alpha_\mathrm{swap}(\teb,\sigma))\delta_{\teb}(B)\pib(\teb)\mub_\mathrm{pr}(\mathrm{d}\teb)\\
	&=\sum_{\sigma\in S_K}\int_{B} r(\teb,\sigma)(1-\alpha_\mathrm{swap}(\teb,\sigma))\delta_{\teb}(A)\pib(\teb)\mub_\mathrm{pr}(\mathrm{d}\teb).
	\end{align}

\end{proof}
\end{proposition}

This generic form of the swapping kernel provides the foundation for both PT and GPT.  We describe these algorithms in the following subsections.

\subsection{The Parallel Tempering case}\label{ss:pt}
We first show how a PT algorithm that only swaps states between the $i^\mathrm{th}$ and $j^\mathrm{th}$ components of the chain can be cast in the general framework presented above. 
To that end, let  $\sigma_{i, j}$ be the permutation of $(1,2,\dots,K),$  which only permutes the $i^\mathrm{th}$ and $j^\mathrm{th}$ components, while leaving the other components invariant (i.e., such that $\sigma(i)=j$, $\sigma(j)=i,$ and $\sigma(k)=k$, $k\neq i,k\neq j$).  We can take $S_K=\{\sigma_{i, j},\ i,j=1,\dots,K\}$ and 
define the PT swapping ratio between components $i$ and $j$ by  $r_{i,j}^\mathrm{(PT)}:\Theta^K\times S_K\mapsto[0,1]$ as 
\begin{align}\label{Eq:pt_swapping_ratio}
r_{i,j}^\mathrm{(PT)}(\teb,\sigma):=\begin{cases}
1 & \text{ if }\sigma=\sigma_{i,j},\\
0 & \text{ otherwise}.
\end{cases}
\end{align}
 Notice that this implies that $r_{i,j}^\mathrm{(PT)}(\teb_\sigma,\sigma^{-1})=r_{i,j}^\mathrm{(PT)}(\teb,\sigma)$ since $\sigma^{-1}_{i,j}=\sigma_{i,j}$ and $r_{i,j}^\mathrm{(PT)}$ does not depend on $\teb$, which in turn  leads to the swapping acceptance probability $\alphaq^\mathrm{(PT)}:\Theta^K\times S_K\mapsto[0,1]$ defined as: 
\begin{align}
\alphaq^\mathrm{{(PT)}}(\teb,\sigma_{i, j})&:=\min\left\{ 1, \frac{\pib(\teb_{\sigma_{i, j}})}{\pib(\teb)}  \right\}, \\ \alphaq^\mathrm{(PT)}(\teb,\sigma)&=0, \ \sigma\neq\sigma_{i,j}.
\end{align}

Thus, we can define the  swapping kernel for the Parallel Tempering algorithm that swaps components $i$ and $j$ as follows:

\begin{definition}[Pairwise Parallel Tempering swapping kernel]	\label{def:pt_ker}
	Let $\teb\in \Theta^K$, $\sigma_{i, j} \in S_K$. We define the \emph{Parallel Tempering swapping kernel}, which proposes to swap states between the $i^\mathrm{th}$ and $j^\mathrm{th}$ chains as $\mathbf{q}_{i, j}^\mathrm{{(PT)}}:\Theta^K\times \mathcal{B}^K\mapsto[0,1] $ given by
	\begin{align}
	&\mathbf{q}_{i, j}^\mathrm{{(PT)}} (\teb, B)=\sum_{\sigma \in S_K}r_{i,j}^\mathrm{{(PT)}}(\teb,\sigma)\left((1-\alphaqt^\mathrm{{{(PT)}}}(\teb,\sigma))\delta_{\teb} (B)\right.\\&\left.+\alphaqt^\mathrm{{{(PT)}}}(\teb,\sigma)\delta_{\teb_\sigma}(B)  \right)\\
	&=\left(1-\min\left\{1,\frac{\pib (\teb_{\sigma_{i,j}})}{\pib(\teb)}   \right\}\delta_{\teb}(B) \right)\\&+\min\left\{ 1,\frac{\pib (\teb_{\sigma_{i,j}})}{\pib(\teb)} \right\}\delta_{\teb_{\sigma_{i,j}}}(B), \quad \forall B\in \mathcal{B}^K.
	\end{align}
\end{definition}

In practice, however, the PT algorithm considers various sequential swaps between chains, which can be understood by applying  the composition of kernels $\mathbf{q}^\mathrm{(PT)}_{i, j}\mathbf{q}^\mathrm{(PT)}_{k, \ell}\dots$ at every swapping step. In its most common form \cite{brooks2011handbook,earl2005parallel,miasojedow2013adaptive}, the PT algorithm, hereafter referred to as standard PT (which on a slight abuse of notation we will denote by PT), proposes to swap states between chains at two consecutive temperatures. Its swapping kernel $\mathbf{q}^\mathrm{(PT)}: \Theta^K\times \mathcal{B}^K\mapsto[0,1]$ is given by 
\begin{align}
\mathbf{q}^\mathrm{(PT)}:=\mathbf{q}^\mathrm{(PT)}_{1, 2}\mathbf{q}^\mathrm{(PT)}_{2, 3} ... \mathbf{q}^\mathrm{(PT)}_{K-1, K}.
\end{align}
Moreover, the algorithm described in \cite{earl2005parallel}, proposes to  swap states every $N_s\geq 1$ steps of MCMC.  The complete kernel for the PT kernel is then given by  \cite{brooks2011handbook,earl2005parallel,miasojedow2013adaptive}  
\begin{align}\label{Eq:vanillapt}
\mathbf{p}^\mathrm{(PT)}:=\mathbf{q}^\mathrm{(PT)}_{1, 2}\mathbf{q}^\mathrm{(PT)}_{2, 3} ... \mathbf{q}^\mathrm{(PT)}_{K-1, K}\mathbf{p}^{N_s},
\end{align}
where $\mathbf{p}$ is a standard reversible Markov transition kernel used to evolve the individual chains independently. 

\rjp{
\begin{remark} 

Although the  kernel $\mathbf{p}$ as well as each of the $\mathbf{q}_{i,i+1}$ are $\mub$-reversible, notice that (\ref{Eq:vanillapt})  does not have a palindromic structure, and as such it is not necessarily $\mub$-reversible. One way of making the PT algorithm reversible with respect to $\bm \mu$  is to consider the palindromic form 
\begin{align}
&\mathbf{p}^\mathrm{(RPT)}:=\\
&\left(\mathbf{q}^\mathrm{(PT)}_{1, 2}\mathbf{q}^\mathrm{(PT)}_{2, 3} ... \mathbf{q}^\mathrm{(PT)}_{K-1, K}\right)\mathbf{p}^{N_s}\left(\mathbf{q}^\mathrm{(PT)}_{K, K-1}...\mathbf{q}^\mathrm{(PT)}_{3, 2}\mathbf{q}^\mathrm{(PT)}_{2, 1}\right),
\end{align}
where RPT stands for \textit{Reversible Parallel Tempering.} In practice, there is not much difference between $\mathbf{p}^\mathrm{(RPT)}$ and $\mathbf{p}^\mathrm{(PT)}$, however, under the additional assumption of geometric ergodicity of the chain (c.f Section \ref{S:AsympInfSwap}) having a reversible kernel is useful to  compute explicit error bounds on the non-asymptotic mean square error of an ergodic estimator \cite{rudolf2011explicit}.
\end{remark}
}

\subsection{Unweighted Generalized Parallel Tempering}\label{SS:uwis}
The idea behind the Unweighted Generalized Parallel Tempering algorithm is to generalize PT so that (i) $N_s=1$ provides a proper mixing of the chains, (ii) the algorithm is reversible with respect to $\mub$, and (iii) the algorithm considers arbitrary sets $S_K$ of swaps (always closed w.r.t inversion), instead of only pairwise swaps. We begin by constructing a kernel of the form  (\ref{Eq:swapping_kernel}). Let $r^\mathrm{{(UW)}}:\Theta^K\times S_K\mapsto [0,1]$ be  a function defined as
\begin{align}\label{Eq:uw_swapping_kernel}
&r^\mathrm{{(UW)}}(\teb,\sigma) := \frac{\pib(\teb_\sigma)}{\sum_{\sigma' \in S_K}\pib(\teb_{\sigma'})}, \quad \teb \in \Theta^K,\sigma\in S_K.
\end{align}  
Clearly, (\ref{Eq:uw_swapping_kernel}) is a swapping ratio according to Definition \ref{def:s_r}. As such,  given some state $\teb\in \Theta^K$,  $r^\mathrm{{(UW)}}(\teb,\sigma)$ assigns a state-dependent probability to each of the ${|S_K|}$ possible permutations in $S_K$. A permutation $\sigma\in S_K$ is then accepted with probability $\alphaq^\mathrm{{(UW)}}(\teb,\sigma)$, given by 
\begin{equation}\label{Eq:uwIS_swapping_acceptance_probability}
\alphaq^\mathrm{{(UW)}}(\teb,\sigma):=\min\left\{ 1,\frac{\pib(\teb_\sigma )r^\mathrm{{(UW)}}(\teb_\sigma,\sigma^{-1})}{\pib(\teb)r^\mathrm{{(UW)}}(\teb,\sigma)} \right\}. 
\end{equation}

Thus, we can define the swapping kernel for the UGPT algorithm, which takes the form of (\ref{Eq:swapping_kernel}), with the particular choice of $r(\teb,\sigma)=r^\mathrm{{(UW)}}(\teb,\sigma)$ and $$\alphaq(\teb,\sigma)=\alphaq^\mathrm{{(UW)}}(\teb,\sigma).$$ Notice that  $\alpha^\mathrm{{(UW)}}_\text{swap}(\teb,\sigma)=1, \forall \sigma\in S_K$.   Indeed, if we further examine Equation \eqref{Eq:uwIS_swapping_acceptance_probability}, we see that 
\begin{align} \frac{\pib(\teb_{\sigma}) r^\mathrm{{(UW)}}(\teb_\sigma,\sigma^{-1})}{\pib(\teb) r^\mathrm{{(UW)}}(\teb,\sigma)} &= \frac{\pib(\teb_\sigma)}{\pib(\teb)}\cdot \frac{\pib(\teb)}{\pib(\teb_\sigma)}\cdot\frac{\sum_{\sigma'}\pib(\teb_{\sigma'})}{\sum_{\hat \sigma}\pib(\teb_{\hat \sigma})} \\&=\frac{\pib(\teb_\sigma)}{\pib(\teb)}\cdot \frac{\pib(\teb)}{\pib(\teb_\sigma)} = 1.\end{align} In practice, this means that the proposed permuted state is always accepted with probability 1. The expression of the UGPT kernel then simplifies as follows. 
\begin{definition}[unweighted swapping kernel]	
	The \emph{unweighted swapping kernel} $\mathbf{q}^\mathrm{{(UW)}}:\Theta^K\times \mathcal{B}^K\mapsto[0,1] $ is defined as
	\begin{equation}
	\mathbf{q}^\mathrm{{(UW)}} (\teb, B)=\sum_{\sigma \in S_K}r^\mathrm{{(UW)}}(\teb,\sigma)\delta_{\teb_\sigma}(B), \ 
	\end{equation}
\end{definition}
$ \forall \teb\in \Theta^K,\  B\in \mathcal{B}^K. $ Applying this swapping kernel successively with the kernel $\mathbf{p}=p_1\times p_2\times \dots p_K $ in the order $\mathbf{q}^\mathrm{{(UW)}} \mathbf{p} \mathbf{q}^\mathrm{{(UW)}}=:\mathbf{p}^\mathrm{{(UW)}}$ gives what we call \emph{Unweighted Generalized Parallel Tempering kernel} $\mathbf{p}^\mathrm{{(UW)}}$. Lastly, we write the UGPT  in operator form as \begin{equation}\label{Eq:uwIS_op}
\mathbf{P}^\mathrm{{(UW)}}:=\mathbf{Q}^\mathrm{{(UW)}} \mathbf{P} \mathbf{Q}^\mathrm{{(UW)}},
\end{equation}
where $\mathbf{P}$ and  $\mathbf{Q}^\mathrm{{(UW)}}$ are the  Markov operators corresponding to the kernels $\mathbf{p}$ and $\mathbf{q}^\mathrm{{(UW)}}$, respectively.  We now investigate the reversibility  of the UGPT kernel. We start with a rather straightforward result.

\begin{proposition}\label{Eq:prop_rev}
	Suppose  that, for any $k=1,2,\dots, K,$ $p_k$ is $\mu_k$-reversible. Then, $\mathbf{p}=p_1\times\dots\times p_K$  is reversible with respect to $\mub$.
	\begin{proof}We prove reversibility by confirming that equation (\ref{Eq:Rev_1}) holds true. To that end, let $\teb\in \Theta^K, \ A,B\in \mathcal{B}^K$, where $A$ and $B$ tensorize, i.e., $A := \prod_{k=1}^K A_k$ and $B := \prod_{k=1}^K B_k$, with $A_1,\ldots, A_K, B_1,\ldots,B_K \in \mathcal{B}(\Theta)$. Then,
		\begin{align} 
		\int_A \pib(\teb) \mathbf{p}(\teb,B)\mathrm{d}\teb&=\prod_{k=1}^K\int_{A_k} \pi(\te_k)p(\te_k,B_k)\mathrm{d}\te_k\\
		&=\prod_{k=1}^K\int_{B_k}\pi(\te_k)p(\te_k,A_k)\mathrm{d}\te_k\\ &=\int_B \pib(\teb) \mathbf{p}(\teb,A)\mathrm{d}\teb. 
		\end{align}
		Showing that the previous equality holds for sets $A, B$ that tensorize is indeed sufficient to show that the claim holds for any $A, B \in \mathcal{B}^K$. This follows from Carath\'eodory's Extension Theorem applied as in the proof of uniqueness of product measures; see \cite[\S 1.3.10, 2.6.3]{Ash2000}, for details.
	\end{proof}
\end{proposition}

We can now prove the reversibility of the chain generated by  $\mathbf{p} ^\mathrm{(UW)}$.

\begin{proposition}[Reversibility of the UGPT chain]\label{proposition:rev_u} Suppose that, for any $k=1,2,\dots, K,$ $p_k$ is $\mu_k$-reversible. Then, the Markov chain generated by  $\mathbf{p} ^\emph{(UW)}$ is \emph{$\mub$}-reversible.
	
\begin{proof}
It follows from Proposition \ref{Prop_MK_Swap} and \ref{Eq:prop_rev} that the kernels $\mathbf{q}^\mathrm{{(UW)}}$ and $\mathbf{p}$ are $\mub$-reversible. Furthermore, since   $\mathbf{p}^\mathrm{{(UW)}}$ is a \textit{palindromic} composition of kernels, each of which is reversible with respect to $\bm \mu$, then, $\mathbf{p}^\mathrm{{(UW)}}$ is reversible with respect to $\bm \mu$ \cite{brooks2011handbook}.

\end{proof}

\end{proposition}

The UGPT  algorithm  proceeds by iteratively  applying the kernel $\mathbf{p}^\mathrm{{(UW)}}$ to a predefined initial state. In particular, states are updated using the procedure outlined in Algorithm \ref{alg:uwIS1}.

%
%
%
%
%

	\begin{algorithm}[tbh]
	\caption{Unweighted Generalized Parallel Tempering.}\label{alg:uwIS1}
	\begin{algorithmic}
		\Function{Generalized Parallel Tempering}{$\mathbf{p},N,\bm \nu $}
		\State Sample $\teb^{(1)}\sim \bm \nu $
		\For {$n=1,2,\dots,N-1$}
		\State \color{gray} \# First swapping kernel \color{black}
		\State Sample $\bm \theta^{(n)}_\sigma \sim \mathbf{q}^\mathrm{{(UW)}}(\teb^{(n)},\cdot)$ 
		\State \color{gray} \# Markov transition kernel $\mathbf{p} $ \color{black}
		\State Sample $\bm {z}^{(n+1)}\sim \mathbf{p} (\bm {\theta}_\sigma^{(n)},\cdot)$ 
		 kernel \color{black}	
		\State \color{gray} \# Second swapping kernel \color{black}		
		\State Sample $\teb^{(n+1)} \sim \mathbf{q}^\mathrm{{(UW)}}(\bm z^{(n+1)},\cdot)$ 
		\EndFor
		\State Output $\{\theta^{(n)}_1\}_{n=1}^{N}$
		\EndFunction
	\end{algorithmic}
	
\end{algorithm}

%
\begin{remark}\label{remark:cost_uw} 
In practice, one does not need to perform $|S_K|$ posterior evaluations when computing $r^\mathrm{(UW)}(\teb^{n},\cdot),$ rather ``just'' $K$ of them. Indeed, since $\pi_j(\te^n_k)\propto\pi(\te_k)^{T_j},$ $k,j=1,2,\dots,K$, we just need to store the values of $\pi(\te^{n}_k),k=1,2,\dots,K$, for a fixed $n,$ and then permute over the temperature indices.

\end{remark}

Let now $\mathcal{Q}:\Theta \mapsto \mathbb{R}$ be a quantity of interest. The posterior mean of $\cal Q$ ,  $\mu(\mathcal{Q}) := \mu_1(\mathcal{Q})$  is approximated using $N \in \mathbb{N}$ samples  by the following ergodic estimator $\widehat{\mathcal{Q}}_{(\mathrm{UW})}$:
$$
\mu(\mathcal{Q})\approx\widehat{\mathcal{Q}}_{(\mathrm{UW})}  = \frac{1}{N-b}\sum_{n=b}^{N} \mathcal{Q}(\te^{(n)}_1) .
$$

\subsubsection{A comment on the pairwise state-dependent PT method of \cite{lkacki2016state}}\label{sss:lacki}

The work \cite{lkacki2016state} presents a similar state-dependent swapping.  We will refer to the method presented therein as Pairwise State Dependent Parallel Tempering (PSDPT). Such a method, however, differs from UGPT  from the fact that (i) only pairwise swaps are considered and (ii) it is not rejection free. We summarize such a method for the sake of completeness. Let $S_{K,\mathrm{pairwise}}$ denote the group of pairwise permutations of $(1,2,\dots,K).$ Given a current state $\teb\in\Theta^K$, the PSDPT algorithm samples a pairwise permutation $\teb_{\sigma_{i, j}}\in S_{K,\mathrm{pairwise}}$ with probability $r_{i,j}^\mathrm{(PSDPT)}(\teb,\sigma_{i,j})$ given by  \begin{align}\label{Eq:pt_swapping_ratio_lacki}
r_{i,j}^\mathrm{(PSDPT)}(\teb,\sigma_{i,j}):= \frac{\exp(-|\Phi(\te_i,y)-\Phi(\te_j;y)|)}{\sum_{k,l}\exp(-|\Phi(\te_k,y)-\Phi(\te_l;y)|)},
\end{align}
and then accepts this swap with probability $$\alpha_\mathrm{swap}^\mathrm{(PSDPT)}(\teb,\sigma_{i,j}):=\min\left\{1,\left(\frac{\pi_1(\te_i)}{\pi_1(\te_j)}\right)^{\frac{1}{T_j}-\frac{1}{T_i}}  \right\}.$$ This method is attractive from an implementation point of view in the sense that it promotes pairwise swaps that have a similar \textit{energy}, and as such, are \textit{likely} (yet not guaranteed) to get accepted. In contrast, UGPT \textit{always} accepts the new proposed state, which in turn leads to a larger amount of  \textit{global} moves, thus providing a more efficient algorithm. This is verified on the numerical experiments.  

\subsection{Weighted Generalized Parallel Tempering}\label{SS:wis}

Following the intuition of the continuous-time Infinite Swapping approach of \cite{dupuis2012infinite,plattner2011infinite}, we propose a second discrete-time algorithm, which we will refer to as \emph{Weighted Generalized Parallel Tempering} (WGPT).  The idea behind this method is to swap the dynamics of the process, that is, the Markov kernels and temperatures, instead of swapping the states such that any given swap is accepted with probability 1.
We will see that the Markov kernel  obtained when swapping the dynamics is not invariant with respect to the product measure of interest $\mub$; therefore, an importance sampling step is needed when computing posterior expectations. 

For a given permutation $\sigma \in S_K$, we define the \textit{swapped Markov kernel} $\mathbf{p}_\sigma : \Theta^K\times \mathcal{B}^K\mapsto[0,1]$ and the \textit{swapped product posterior measure} $\mub_\sigma$ (on the measurable space $(\Theta^K,\mathcal{B}^K)$) as: 
\begin{align}
&\mathbf{p}_\sigma(\teb,\cdot) = {p}_{\sigma(1)}(\theta_1,\cdot) \times \dots\times {p}_{\sigma(K)}(\theta_K,\cdot),\\
&\bm \mub_\sigma:=\mu_{\sigma(1)}\times\dots\times\mu_{\sigma(K)},
\end{align}
where the swapped posterior measure has a density with respect to $\bm \mu_\textrm{prior}$ given by 
\begin{align} \label{weighs_wIS}
&\bm \pib_\sigma(\teb):=\pi_{\sigma(1)}(\te_1)\times\dots\times\pi_{\sigma(K)}(\te_K), \quad \te\in \Theta^K,\sigma\in S_K
\end{align}
Moreover, we define the swapping weights
 \begin{align}\label{Eq:weighted_is_weights}
	w_\sigma(\teb):=\frac{\pib_\sigma(\teb)}{\sum_{\sigma'\in S_K}\pib_{\sigma'}(\teb)}, \quad \teb\in \Theta^K, \ \sigma \in S_K.
\end{align}
Note that, in general, $\pib_\sigma(\teb)\neq\pib(\teb_\sigma)$, and as such, $w_\sigma(\teb)\neq r^\mathrm{{(UW)}}(\teb,\sigma)$, with $w_\sigma$ defined as in \eqref{Eq:weighted_is_weights}.

\begin{definition}
We define the \emph{Weighted Generalized Parallel Tempering} kernel $\mathbf{p}^\mathrm{(W)} : \Theta^K\times B^K\mapsto[0,1]$ as the following state-dependent, convex combination of kernels: 

\begin{equation} \mathbf{p}^{\mathrm{(W)}}(\bm \te, \cdot) := \sum_{\sigma \in S_K} w_\sigma(\teb) \mathbf{p}_\sigma(\bm \te,\cdot), \quad \teb\in \Theta^K, \ \sigma\in S_K.\end{equation} 
\end{definition}
Thus, the WGPT chain is obtained by iteratively applying $\mathbf{p}^{\mathrm{(W)}}$.
 We show in proposition \ref{proposition:rev} that the resulting Markov chain has invariant measure  
 \begin{align}
 \bm \mub _\text{W}&=\frac{1}{{|S_K|}}\sum_{\sigma\in S_K} \bm \mub_\sigma=\tilde \mu \times \dots \times \tilde \mu, \label{EQ:meassure_is}
 \end{align}
 with $\tilde \mu=\frac{1}{|S_K|}\sum_\sigma \mu_\sigma,$ i.e., the average with tensorization. Furthermore, $\mubis$ has a density (w.r.t the prior $\mub^0$) given by
 \begin{align}
\bm \pib_\text{W}(\bm \te) &=\frac{1}{{|S_K|}}\sum_{\sigma\in S_K} \pib_\sigma(\bm \te), \quad \teb\in \Theta^K,
\end{align} and a similar average and then tensorization representation applies to $\pib_{\mathrm{W}}.$ We now proceed to show that  $\mathbf{p}^{(\mathrm{W})}(\teb,\cdot)$ is $\mubis$-reversible (hence $\mubis$-invariant).


\begin{proposition}[Reversibility of the WGPT chain]\label{proposition:rev} Suppose that, for any $k=1,2,\dots, K$ $p_k$ is $\mu_k$-reversible. Then, the Markov chain generated by  $\mathbf{p} ^\emph{(W)}$ is \emph{$\mub_\textrm{W}$}-reversible.
	%
	%
	\begin{proof}
		 We show reversibility by showing that (\ref{Eq:Rev_1}) holds true. Thus, for $\teb\in \Theta^K,\  A,B\in \mathcal{B}^K$, with $A:=A_1\times\dots\times A_K$, $A_k\in \mathcal{B}(\Theta)$, and with $B_k$ defined in a similar way, we have that: 		
		\begin{align}
		&\int_A \mathbf{p}^{\mathrm{(W)}}(\teb,B)\bm \pi _\text{W}(\teb)\mub_\mathrm{pr}(\mathrm{d}\teb)\\	&=\int_A\left[\sum_{\sigma \in S_K} w_\sigma(\teb)\mathbf{p}_\sigma(\bm \te, B)\right]  \frac{\sum_{\rho \in S_K}\pib_\rho(\teb)}{{|S_K|}}\mub_\mathrm{pr}(\mathrm{d}\teb)\\
		&=\int_A\left[\sum_{\sigma \in S_K} \frac{\bm \pi_\sigma(\bm \te)}{\sum_{\sigma' \in S_K}\pib_{\sigma'}(\teb)}\mathbf{p}_\sigma(\bm \te, B)\right] \\& \times\frac{\sum_{\rho \in S_K}\pib_\rho(\teb)}{{|S_K|}}\mub_\mathrm{pr}(\mathrm{d}\teb)\\
		&=\frac{1}{{|S_K|}}\sum_{\sigma\in S_K}\int_A  \pib_\sigma(\teb)\mathbf{p}_\sigma(\bm \te, B)\mub_\mathrm{pr}(\mathrm{d}\teb)=I.
		\end{align}
		\rjp{From proposition \ref{Eq:prop_rev}, and multiplying and dividing by $$\sum_{\rho \in S_K}\pib_{\rho}(\teb)$$  we obtain}
		\begin{align}
		&I=\frac{1}{{|S_K|}}\sum_{\sigma\in S_K}\int_B  \pib_\sigma(\teb)\mathbf{p}_\sigma(\bm \te, A)\mub_\mathrm{pr}(\mathrm{d}\teb)
		\quad \text{(by Prop. \ref{Eq:prop_rev})}\\
		&=\frac{1}{{|S_K|}}\sum_{\sigma\in S_K}\int_B  \frac{\pib_\sigma(\teb)\mathbf{p}_\sigma(\teb,A)}{\sum_{\sigma' \in S_K}\pib_{\sigma'}(\teb)}\sum_{\rho \in S_K}\pib_{\rho}(\teb)\mub_\mathrm{pr}(\mathrm{d}\teb)
		\\
		&=\sum_{\sigma\in S_K}\int_B  w_\sigma(\teb)\mathbf{p}_\sigma(\bm \te, A)\pib_\mathrm{W}(\teb)\mub_\mathrm{pr}(\mathrm{d}\teb)\\
		&=\int_B \mathbf{p}^{\mathrm{(W)}}(\teb,A)\bm \pi _\text{W}( \teb)\mub_\mathrm{pr}(\mathrm{d}\teb).
		\end{align}
		where once again, in light of Carath\'eodory's Extension Theorem, it is sufficient to show that reversibility holds for sets that tensorize. 
	\end{proof}
\end{proposition}

 We remark that the measure  $\mubis$ is not of interest per se.
However, we can use importance sampling to compute posterior expectations. Let $\mathcal{Q}(\teb):=\mathcal{Q}(\te_1)$ be a $\mub$-integrable quantity of interest. We can write 
 \begin{align}
\E_{\mu_1}[\mathcal{Q}]&=\E_{\bm \mu}[\mathcal{Q}(\te_1)]=\E_{\bm \mu_\text{W}}\left[\mathcal{Q}(\te_1)\frac{\pib(\teb)}{\pib_\text{W}(\teb)}\right] \\&=\frac{1}{{|S_K|}}\sum_{\sigma\in S_K}\E_{\bm \mu_\text{W}}\left[\mathcal{Q}(\te_{\sigma(1)})\frac{\pib(\teb_\sigma)}{\pib_\text{W}(\teb_\sigma)}\right]. \label{Eq:imp_sampling}
\end{align}

The last equality can be justified since $\mubis$ is invariant by permutation of coordinates. Thus, we  can define the following (weighted) ergodic estimator $\widehat{\mathcal{Q}}_{(\mathrm{W})} $ of the posterior mean of a quantity of interest $\mathcal{Q} $ by 
\begin{align}
\mu(\mathcal{Q})&\approx\\&\widehat{\mathcal{Q}}_{(\mathrm{W})} = \frac{1}{{|S_K|}}\frac{1}{N} \sum_{\sigma\in S_K} \sum_{n=1}^N \frac{\pib(\teb^{(n)}_\sigma)}{\pib_\text{W}(\teb^{(n)}_\sigma)}\mathcal{Q}(\te^{(n)}_{\sigma(1)})\\
&= \frac{1}{{|S_K|}}\frac{1}{N} \sum_{\sigma\in S_K} \sum_{n=1}^N \widehat{w}(\teb^{(n)},\sigma)\mathcal{Q}(\te^{(n)}_{\sigma(1)}),\label{Eq:weighted_samples}
\end{align}
where we have denoted the importance sampling weights by $\widehat{w}(\teb,\sigma):=\frac{\pib(\teb_\sigma)}{\pib_\text{W}(\teb_\sigma)}=\frac{\mathrm{d}\mub}{\mathrm{d}\mubis}(\teb_\sigma)$ and where $N$ is the number of samples in the chain. Notice that $w(\teb,\sigma)=\widehat{w}(\teb,\sigma^{-1}).$ As a result, the WGPT algorithm produces an estimator based on $N K$  weighted samples, rather than ``just'' $N$, at the same computational cost of UGPT. Thus, the previous estimator evaluates the quantity of interest $\mathcal{Q}$ not only in the points $\mathcal{Q}(\te_1^{(n)})$, but also in all states of the parallel chains, $\mathcal{Q}(\te^{(n)}_{\sigma(1)})$ for all $\sigma\in S_K$, namely $\mathcal{Q}(\te^{(n)}_k), \ k=1,2,\dots,K$. 
\rjp{
\begin{remark}
Although it is known that, in some cases, an importance sampling estimator can be negatively affected by the dimensionality of the parameter space $\Theta$ (see e.g., \cite[Remark 1.17]{asmussen2007stochastic} or \cite[Examples 9.1-9.3]{mcbook}), we argue that this is not the case for our estimator. Indeed, notice that the importance-sampling weights $\widehat{w}(\teb,\sigma)$  are always upper bounded by $|S_K|$, and do not blow up when the dimension goes to infinity. In Section \ref{ss:high_dim} we present a numerical example on a high-dimensional problem. The results on that section evidence the robustness of WGPT with respect to the dimension of $\te$.  
\end{remark}
}

The Weighted Generalized Parallel Tempering procedure is shown in Algorithm \ref{alg:wIS}. To reiterate, we remark that sampling from $\mathbf{p}_\sigma(\teb^{(n)},\cdot)$ involves a swap of dynamics, i.e., kernels and temperatures.

	\begin{algorithm}[tbh]
	\caption{Weighted Generalized Parallel Tempering.}\label{alg:wIS}
	\begin{algorithmic}
		\Function{Weighted Generalized Parallel Tempering}{$\ \mathbf{p},N,\bm \nu $}
		\State Sample $\teb^{(1)}\sim \bm \nu $
		\For {$n=1,2,\dots,N-1$}
		\State \color{gray} \# Sample permutation $\sigma$ with probability $w_{\sigma}(\teb^n)$ \color{black}
		\State Sample $\sigma \sim \{w_{\sigma'}(\teb^n) \}_{\sigma'\in S_K}$ 
		\State \color{gray} \# Sample state with the swapped Markov kernel \color{black}
		\State Sample $\bm {\theta}^{(n+1)}\sim \mathbf{p}_\sigma (\bm {\theta}^{(n)},\cdot)$ 
		\EndFor
		\State Output $\{\teb^{(n)}\}_{n=1}^{N}, \{\{w_{\sigma'}(\teb^n) \}_{\sigma'\in S_K}\}_{n=1}^N$.
		\EndFunction
	\end{algorithmic}
	
\end{algorithm}

Just as in Remark \ref{remark:cost_uw}, one only needs to evaluate the posterior $K$ times (instead of $|S_K|$) to compute $w_{(\cdot)}(\teb^n)$.

\section{Ergodicity of Generalized Parallel Tempering} \label{S:AsympInfSwap}

\subsection{Preliminaries} \label{ss:prelim}
We assume that the chains generated by the MCMC  kernels $p_k$, for $k=1,\dots,K$, are aperiodic,  $\mu_k$-irreducible \cite{asmussen2007stochastic}, and have invariant measure $\mu_k$ on the measurable space $(\Theta,\mathcal{B}(\Theta))$. Let  $r\in[1,\infty)$ and $\mu\in \mathcal{M}(\Theta)$ be a ``reference'' probability measure. On a BIP setting, this reference measure is considered to be the posterior. We define the following spaces
\begin{align}
&L_r=L_r(\Theta,\mu)=\left\{f:\Theta\mapsto\R, \ \mu\text{-measurable, } \right. \\ &\left. \quad \quad \quad \quad  \text { s.t } \lno f\rno_r^r\defeq\int |f(\te)|^r\mu(\mathrm{d}\te)<\infty       \right\},\\
&L_r^0=L^0_r(\Theta,\mu)=\left\{f\in L_r(\Theta,\mu), \right. \\&\left. \quad \quad \quad \quad   \text{ s.t } \mu(f):=\int_\Theta f(\te)\mu(\mathrm{d}\te)=0  \right\}.
\end{align}
Moreover, when $r=\infty$, we define \begin{equation}
L_\infty(\Theta,\mu):=\left\{f:\Theta\mapsto \R, \ \text{ s.t } \inf_{\underset{B\in \mathcal{B}(\Theta)}{\mu(B)=0}}\underset{y\in \Theta\backslash B}{\sup} \abs{f(y)}<\infty \right \}.
\end{equation}

\noindent Notice that, clearly, ${L^0_r(\Theta,\mu)}\subset{L_r(\Theta,\mu)}$. In addition we define the spaces of measures \begin{align}
&\mathcal M_r(\Theta,\mu):=\{ \nu \in \mathcal M(\Theta) \text{ s.t. } \nu\ll\mu, \  \lno \nu\rno_{L_r(\Theta,\mu)}< \infty \},\\
&\text{where}\quad \lno \nu \rno_{L_r(\Theta,\mu)}\defeq \lno\frac{\mathrm{d}\nu}{\mathrm{d}\mu} \rno_{L_r(\Theta,\mu)}.
\end{align}

\noindent Notice that the definition of $L_r$-norm depends on the reference measure $\mu,$ and on $\Theta$.

A Markov operator  $P:L_r(\Theta,\mu)\mapsto L_r(\Theta,\mu)$ with invariant measure $\mu$ is a bounded linear operator whose norm is given by  \begin{equation}
\lno P \rno_{L_r(\Theta,\mu)\mapsto L_r(\Theta,\mu)} :=\sup_{\lno f \rno_{L_r(\Theta,\mu)}=1}  \lno  Pf \rno_{L_r(\Theta,\mu)}, 
\end{equation}
for $\ f\in L_r(\Theta,\mu).$ It is well-known (see, e.g., \cite{rudolf2011explicit}) that any Markov operator $P$ on $L_r(\Theta,\mu)$ with invariant measure $\mu$ can be understood as a weak contraction in $L_r(\Theta,\mu),$  i.e., $\lno P\rno_{L_r(\Theta,\mu)\mapsto L_r(\Theta,\mu)}\leq 1$. To quantify the convergence of a Markov chains generated by a Markov operator $P$, we define the concept of geometric ergodicity.  Let $r\in[1,\infty]$.  A Markov operator $P$ with invariant measure $\mu\in \mathcal{M}(\Theta)$ is said to be \emph{$L_r(\Theta,\mu)$-geometrically ergodic} if for all probability measures $\nu\in \mathcal{M}_r(\Theta,\mu)$ there exists an $\alpha\in(0,1)$ and $C_\nu<\infty$ such that \begin{align}
\lVert \nu P^n-\mu\rVert_{L_r(\Theta,\mu)}\leq C_\nu\alpha^n, \ \ n\in\ \mathbb{N}.\label{eq:l2ge}
\end{align}
A related concept to $L_2$-geometric ergodicity is that of $L_2$-spectral gap. A Markov operator $P: {L_2(\Theta,\mu)}\mapsto L_2(\Theta,\mu)$ with invariant measure $\mu\in \mathcal{M}(\Theta)$  has an \emph{$L_2(\Theta,\mu)$-spectral gap} $1-\beta>0$, with $\beta<1$, if the following holds 
\begin{equation}
\lno P \rno_{L^0_2(\Theta,\mu)\mapsto L^0_2(\Theta,\mu)}\leq\beta. \label{eq:l2sg}
\end{equation}

The next Proposition, whose proof can be found e.g., in \cite{rudolf2011explicit}, relates the existence of an $L_2$-spectral gap to the geometric ergodicity of the chain (with $\beta\leq \alpha$, in general). 

\begin{proposition}\label{proposition:l2sg_rge}
Let $P:L_2(\Theta,\mu)\mapsto L_2(\Theta,\mu)$ be a $\mu$ reversible Markov transition operator. The existence of an $L_2(\Theta,\mu)$-spectral gap implies $L_r(\Theta,\mu)$-geometric ergodicity for any $r\in[1,\infty]$.  
\begin{proof}
	The previous claim is shown in   \cite[Proposition 3.17 and Appendix A.4]{rudolf2011explicit}. It is also shown in \cite{rudolf2011explicit} that, in general, $\beta\leq \alpha$, with $\alpha,\beta$ given as in Equations \eqref{eq:l2ge} and \eqref{eq:l2sg}.
	\end{proof}
\end{proposition}

%
Our path to prove ergodicity of the GPT algorithms will be to show the existence of an $L_2$-spectral gap. 

\subsection{Geometric ergodicity and $L_2$-spectral gap for GPT}\label{ss:main_results}

The main results of this section are presented in Theorem \ref{proposition:uwIS_ergodicity} and Theorem \ref{proposition:wIS_ergodicity}, which show the existence of an $L_2$-spectral gap for both the UGPT and WGPT algorithms, respectively. 

We begin with the definition of overlap between two probability measures. Such a concept will later be used to bound the spectral gap of the GPT  algorithms. 

\begin{definition}[Density overlap]\label{def:overlap}
	Let $\mu_k,\mu_j$ be two probability measures on the measurable space $(\Theta,\mathcal{B}(\Theta))$, each having respective densities $\pi_k(\te),\pi_j(\te), \ \te\in \Theta,$ with respect to some common reference measure $\nu_\Theta$ also on $(\Theta,\mathcal{B}(\Theta))$. We define the \emph{overlap} between $
	\pi_k(\te)$ and $\pi_j(\te)$	as  \begin{align}
	\eta_{\nu_\Theta}(\pi_k,\pi_j):&=\int_\Theta\min \{\pi_k(\te),\pi_j(\te)\}\nu_\Theta(\mathrm{d\te}) \\&=1-\frac{1}{2}\lno \mu_k-\mu_j\rno_{L_1(\Theta,\nu_{\Theta})}.
	\end{align} An analogous definition holds for $\pib_\sigma,\pib_\rho$, with $\rho,\sigma\in S_K$.
\end{definition}

\begin{assumption}
	
\label{All_assumptions}For $k=1,\dots,K$, let $\mu_k\in\mathcal M_1(\Theta,\mu_\mathrm{pr})$ be given as in \eqref{Eq:RN_single},  $p_k:\Theta\times \mathcal{B}(\Theta)\mapsto[0,1]$ be the Markov kernel associated to the $k^\text{th}$ dynamics and let  $P_k:L_r(\Theta,\mu_k)\mapsto L_r(\Theta,\mu_k)$ be its corresponding $\mu_k$ invariant Markov operator. In addition, for $\sigma,\rho\in S_K$, define the measures $\mub_\sigma,\mub_\rho\in\mathcal{M}(\Theta^K)$ as in Equation (\ref{Eq:Product_measure}). Throughout this work it is assumed that:
	\begin{enumerate}[label={{C\arabic*}}.,ref=C\arabic*]
		\item\label{ass:1}The Markov kernel $p_k$ is $\mu_k$-reversible. 
		\item\label{ass:2} The Markov operator $P_k$ has an  $L_2(\Theta,\mu_k)$ spectral gap.
		\item\label{ass:3} For any $\sigma,\rho\in S_K$, $\Lambda_{\sigma,\rho}:=\eta_{\mub_\mathrm{pr}}(\pib_\sigma,\pib_\rho)>0,$ with $\pib_\sigma,\pib_\rho$ defined as in \eqref{weighs_wIS}.
	\end{enumerate}

\end{assumption}

These assumptions are relatively mild. In particular, \ref{ass:1} and 
	\ref{ass:2} are known to hold for many commonly-used Markov transition kernels, such as RWM, Metropolis-adjusted Langevin Algorithm, Hamiltonian Monte Carlo, (generalized) preconditioned Crank-Nicolson, among others, under mild regularity conditions on $\pi$ \cite{asmussen2007stochastic,hairer2014spectral}.  Assumption \ref{ass:3} holds true given the construction of the product measures in Section \ref{S:pt_e_is}.  

 We now present an auxiliary result that we will use to bound the spectral gap of both the Weighted and Unweighted GPT algorithms.

\begin{proposition}\label{proposition:geom_erg_joint}
	Suppose that Assumption \ref{All_assumptions}  holds and let $\mathbf{P}:=\bigotimes_{k=1}^KP_k: L_2(\Theta^K,\mub)\mapsto L_2(\Theta^K,\mub),$ with invariant measure $\mub=\mu_1\times\dots\times\mu_K$. Then, $\mathbf{P}$ has an $L_2(\Theta^K,\mub)$-spectral gap, i.e., $\lno \mathbf{P} \rno_{L^0_2(\Theta^K,\bm \mu)\mapsto L^0_2(\Theta^K,\mub)}<1$. Moreover, the Markov chain obtained from  $\mathbf{P}$ is $L_r$ geometrically ergodic, for any $r\in[1,\infty]$. 
	\begin{proof}
		We limit ourselves to the case $K=2$, since the case for $K>2$ follows by induction. Denote by $I:L_2(\Theta,\mu_k)\mapsto L_2(\Theta,\mu_k), k=1,2$ the identity Markov transition operator,  and let $f\in L_2(\Theta^2,\mub).$ Notice that $f$ admits a spectral representation in $L_2(\Theta^2,\mub)$ given by $f(\teb)=\sum_{k,j}\phi_k(\te_1)\psi_j(\te_2)c_{k,j}$, with $c_{k,j}\in \R$, and   where $\{ \phi_k \}_{i\in\mathbb{N}}$ is a complete orthonormal basis (CONB) of $L_2(\Theta,\mu_1)$ and $\{\psi_j\}_{j\in\mathbb{N}}$ is a CONB of $L_2(\Theta,\mu_2)$, so that $\{\phi_k\otimes\psi_j\}_{k,j\in\mathbb{N}}$ is a CONB of $L_2(\Theta^2,\mub)$. Moreover, we assume  that $\phi_0=\psi_0=1$, and write, for notational simplicity  $\lno P_1 \rno=\lno P_1\rno_{L_2(\Theta,\mu_1)\mapsto L_2(\Theta,\mu_1)}$, and $\lno{P_2}\rno=\lno P_2\rno_{L_2(\Theta,\mu_2)\mapsto L_2(\Theta,\mu_2)}$.  Lastly, denote $f_0=f-c_{0,0}$, so that $f_0\in L_2^0(\Theta^2,\mub)$.  Notice that 
	\begin{align}
&\lno (P_1 \otimes I)f_0\rno ^2_{L_2(\Theta^2,\mub)} =\lno \sum_{(k,j)\neq(0,0)}\left(P_1\phi_k\right)\psi_jc_{k,j} \rno^2_{L_2(\Theta^2,\mub)}\\
&=\lno \sum_{j=0}^{\infty}\left(\sum_{k=1}^{\infty}P_1\phi_kc_{k,j}\right)\psi_j+\sum_{j=1}^\infty c_{0,j}P_1\phi_0\psi_j\rno^2_{L_2(\Theta^2,\mub)}.\label{eq:p1}
\end{align}
\rjp{Splitting the sum, we get from the orthonormality of the basis that:}
\begin{align}
	&\eqref{eq:p1}=\sum_{j=1}^\infty \lno \sum_{k=1}^\infty P_1\phi_kc_{k,j}+c_{0,j}P_1\phi_0\rno^2_{L_2(\Theta,\mu_1)}\\ &+\lno \sum_{k=1}^\infty P_1\phi_kc_{k,0}\rno^2_{L_2(\Theta,\mu_1)} \\
	&=\sum_{j=1}^\infty \lno P_1\left(\sum_{k=1}^\infty \phi_kc_{k,j}\right)\rno ^2_{L_2(\Theta,\mu_1)} +\sum_{j=1}^\infty\lno c_{0,j}\phi_0\rno^2_{L_2(\Theta,\mu_1)}\\&+\lno P_1\left(\sum_{k=1}^\infty \phi_kc_{k,0}\right)\rno^2_{L_2(\Theta,\mu_1)}\\
	&\leq \sum_{j=1}^{\infty}\left(\lno P_1\rno ^2  \sum_{k=1}^\infty c_{k,j}^2+c_{0,j}^2\right)+\lno P_1\rno^2  \sum_{k=1}^{\infty}c_{k,0}^2\\
	&=\lno P_1\rno ^2  \lno f_0\rno_{L_2(\Theta^2,\mub)}^2 +(1-\lno P_1\rno  ^2)\sum_{j=1}^\infty (c_{0,j})^2.  
	\end{align} 
Proceeding similarly, we can obtain an equivalent bound for $\lno (I\otimes P_2)f_0\rno_{L_2(\Theta^2,\mub)} ^2$. We are now ready to bound $\lno \mathbf{P}\rno ^2_{L_2(\Theta^2,\mub)\mapsto L_2(\Theta^2,\mub)}$:
\begin{align}
&\lno \mathbf{P}f_0\rno ^2_{L_2(\Theta^2,\mub)}= \lno (P_1\otimes P_2) f_0\rno_{L_2(\Theta^2,\mub)} ^2\\&=\lno (P_1\otimes I)(I\otimes P_2)f_0\rno_{L_2(\Theta^2,\mub)}^2\\
&\leq \lno P_1\rno  ^2 \lno (I\otimes P_2) f_0\rno_{L_2(\Theta^2,\mub)}^2\\ &+(1-\lno P_1\rno  ^2 )\\&\times\left(\sum_{j=1}^\infty \left((I\otimes P_2)\sum_{(\ell,k)\neq (0,0)}c_{\ell,k}\phi_\ell\psi_k,\phi_0\psi_j\right)^2\right)\\
&= \lno P_1\rno  ^2 \lno (I\otimes P_2) f_0\rno_{L_2(\Theta^2,\mub)}^2\\&+(1-\lno P_1\rno  ^2 )\left(\sum_{j=1}^\infty \left(\sum_{k=1}^\infty c_{0,k}(P_2\psi_k),\psi_j \right)^2\right)\\
&\leq \lno P_1\rno  ^2 \lno (I\otimes P_2) f_0\rno_{L_2(\Theta^2,\mub)}^2\\&+(1-\lno P_1\rno  ^2 )\lno P_2\left(\sum^\infty_{k=1}c_{0,k}\psi_k\right)\rno_{L_2(\Theta,\mu_2)}^2
\end{align}
\begin{align}
&\leq \lno P_1\rno ^2 \lno P_2\rno ^2 \lno f_0\rno_{L_2(\Theta^2,\mub)}^2\\
&+\lno P_1\rno  ^2 (1-\lno P_2\rno  ^2)\left(\sum_{j=1}^\infty c_{j,0}^2\right)\\
	&+(1-\lno P_1\rno  ^2)\lno P_2 \rno ^2 \left(\sum_{k=1}^\infty c_{0,k}^2\right)
\end{align}
Assuming without loss of generality that $\lno P_1\rno  \geq \lno P_2\rno $, we can use the inequality above to bound
\begin{align}
\lno \mathbf{P}f_0\rno ^2_{L_2(\Theta^2,\mub)}&\leq \lno P_1\rno ^2 \lno P_2\rno ^2 \lno f_0\rno_{L_2(\Theta^2,\mub)}^2 \\ &+ \lno P_1\rno ^2 (1-\lno P_2\rno  ^2)\underbrace{\left(\sum_{j=1}^\infty c_{j,0}^2+\sum_{k=1}^\infty c_{0,k}^2\right)}_\text{$\leq\lno f_0\rno_{L_2(\Theta^2,\mub)} ^2$}\\&\leq \lno P_1\rno ^2 \lno f_0 \rno_{L_2(\Theta^2,\mub)}^2.
\end{align}		
Thus, we have that $$\lno \bm P \rno_{L^0_2(\Theta^2,\mub)\mapsto L^0_2(\Theta^2,\mub)} \leq \max_{k=1,2}\{\lno P_k\rno_{L^0_2(\Theta,\mu_k)\mapsto L^0_2(\Theta,\mu_k)} \}<1.$$ The previous result can easily be extended to $K>2$. Lastly, $L_r(\Theta^K,\mub)$-geometric ergodicity $\forall r\in[1,\infty]$ follows from proposition \ref{proposition:l2sg_rge}.
	\end{proof}
\end{proposition}

We can use the previous result to prove the geometric ergodicity of the UGPT  algorithm:

\begin{thm}[Ergodicity of UGPT ]\label{proposition:uwIS_ergodicity} Suppose \emph{Assumption  \ref{All_assumptions}} holds and denote by $\mub$ the invariant measure of the UGPT Markov operator $\mathbf{P}^{\emph{(UW)}}$. Then, $\mathbf{P}^{\emph{(UW)}}$ has an $L_2(\Theta^K,\mub)$-spectral gap. Moreover, the chain generated by $\mathbf{P}^{\emph{(UW)}}$  is $L_r(\Theta^K,\mub)$-geometrically ergodic for any $r\in[1,\infty]$.

	\begin{proof}
Recall that $\mathbf{P}^\mathrm{(UW)}:=\mathbf{Q}^\mathrm{(UW)} \mathbf{P} \mathbf{Q}^\mathrm{(UW)}$. From the definition of operator norm, we have that
\begin{align}
&\lno \mathbf{P}^{\mathrm{(UW)}}\rno_{L^0_2(\Theta^K,\mub)\mapsto L^0_2(\Theta^K,\mub)}\\&\leq \lno \mathbf{Q}^{\mathrm{(UW)}} \rno^2_{L^0_2(\Theta^K,\mub)\mapsto L^0_2(\Theta^K,\mub)}\lno \mathbf{P} \rno_{L^0_2(\Theta^K,\mub)\mapsto L^0_2(\Theta^K,\mub)}\\
& \leq  \lno \mathbf{P} \rno_{L^0_2(\Theta^K,\mub)\mapsto L^0_2(\Theta^K,\mub)}<1,
\end{align}
where the previous line follows from Proposition \ref{proposition:geom_erg_joint} and the fact that $\mathbf{Q}^{\mathrm{(UW)}}$ is a weak contraction in $L_2(\Theta^K,\mub)$ (see, e.g., \cite[Proposition 1]{baxter1995rates}). Lastly, $L_r(\Theta^K,\mub)$-geometric ergodicity $\forall r\in[1,\infty]$ follows from  proposition \ref{proposition:l2sg_rge} and the fact that $\mathbf{P}^\mathrm{(UW)}$ is $\mub$-reversible by Proposition \ref{proposition:rev_u}.
	\end{proof}
\end{thm}
We now turn to proving geometric ergodicity for the WGPT algorithm. We begin with an auxiliary result, lower-bounding the variance of a $\mubis$-integrable functional $f\in L_2(\Theta^K,\mubis)$. 

\begin{proposition}\label{proposition:bound_var}Let $f\in L^0_2(\Theta^K,\mubis)$ be a $\mubis$-integrable function such that $\lno f \rno_{L_2(\Theta^K,\mubis)}=1,$ and denote by $\V_{\mubis}[f]$, $\V_{\mub_\sigma}[f]$ the variance of $f$ with respect to $\mubis,\mub_\sigma$, respectively with $\sigma\in S_K$. In addition, suppose Assumption \ref{All_assumptions} holds. Then, it can be shown that
\begin{align}
0<\frac{\Lambda_m}{2-\Lambda_m}\leq \frac{1}{{|S_K|}}\sum_{\sigma \in S_K}\V_{\mub_\sigma}[f]\leq \V_{\mubis}[f]=1,
 \end{align}
with $\Lambda_m=\underset{\sigma,\rho\in S_K}{\min}\{\Lambda_{\sigma,\rho}\}$ and $\Lambda_{\sigma,\rho}$ as in Assumption \ref{ass:3}.
\begin{proof}
	The proof is technical and tedious and is presented in \hyperref[appn]{Appendix} \ref{AP:proof_bound_cov}.
	\end{proof}
\end{proposition}

We are finally able to prove the ergodicity of the WGPT algorithm. 

\begin{thm}[Ergodicity of WGPT] \label{proposition:wIS_ergodicity}Suppose Assumption \ref{All_assumptions} holds for some $r\in[1,\infty]$ and denote by $\mub_\emph{W}$ the invariant measure of the WGPT  Markov operator $\mathbf{P}^{\emph{(W)}}$. Then, $\mathbf{P}^{\emph{(W)}}$ has an $L_2(\Theta^K,\mub_\emph{W})$-spectral gap. Moreover, the chain generated by $\mathbf{P}^{\emph{(W)}}$ is $L_r(\Theta^K,\mub_\emph{W})$ geometrically ergodic for any $r\in[1,\infty]$. 
	
	\begin{proof}
	Let $\mathcal{L}:= \{f\in L_2^0(\Theta^K,\mubis) :  \lno f\rno_{L_2^0(\Theta^K,\mubis)}=1\}$, and, for notational clarity, write $$\lno \mathbf{P}_\sigma \rno_{L^0_2}:=\lno \mathbf{P}_\sigma^2 \rno_{L^0_2(\Theta^K,\bm \mub_\sigma)\mapsto L^0_2(\Theta^K,\bm \mub_\sigma)}.$$ Then, from the definition of operator norm,
	\begin{align}
	&\lno \mathbf{\Pw} \rno_{L^0_2(\Theta^K,\bm \mubis)\mapsto L^0_2(\Theta^K,\bm \mubis) }^2\\ &={\sup}_{f \in \mathcal{L}} \lno \mathbf{\Pw} f \rno^2_{L_2(\Theta^K,\bm \mubis)}	\\
	&=\underset{f \in \mathcal{L}}{\sup}\int_{\Theta^K} \abs{\sum_{\sigma \in S_K} w_\sigma (\teb) \int_{\Theta^K} f(\bm y)\mathbf{p}_\sigma(\teb,\mathrm{d}\bm y) }^2 \mubis(\mathrm{d}\teb)\\	
	&\leq \underset{f \in \mathcal{L}}{\sup} \int_{\Theta^K} \sum_{\sigma\in S_K} w_\sigma(\teb)\abs{\int_{\Theta^K}f(\bm y) \mathbf{p}_\sigma(\teb,\mathrm{d}\bm y)}^2\mubis(\mathrm{d}\teb)  \label{Eq:from_convexity}\\
	&=\underset{f \in \mathcal{L}}{\sup} \frac{1}{{|S_K|}} \sum_{\sigma \in S_K} \int_{\Theta^K} \abs{ \int_{\Theta^K} f(\bm y)\mathbf{p}_\sigma(\teb,\mathrm{d}\bm y) }^2 \mub_\sigma(\mathrm{d} \teb), \label{Eq:ine_sg_weighted}
	\end{align}
	where the second to last line follows from the convexity of $(\cdot)^2$ and the last line follows from the definition of $w_\sigma$ and $\mubis.$ Now, let $\fbs:=\mub_\sigma(f)$.  Notice that we have 
	\begin{align}
	&\int_{\Theta^K} \abs{ \int_{\Theta^K} f(\bm y)\mathbf{p}_\sigma(\teb,\mathrm{d}\bm y) }^2 \mub_\sigma(\mathrm{d} \teb)\\
	=&\int_{\Theta^K} \abs{ \int_{\Theta^K} (f(\bm y)-\fbs +\fbs)\mathbf{p}_\sigma(\teb,\mathrm{d}\bm y) }^2 \mub_\sigma(\mathrm{d} \teb)\\
	=&\int_{\Theta^K} \left(\abs{ \int_{\Theta^K} (f(\bm y)-\fbs)\mathbf{p}_\sigma(\teb,\mathrm{d}\bm y)}^2+\abs{\int_{\Theta^K}\fbs\mathbf{p}_\sigma(\teb,\mathrm{d}\bm y) }^2\right. \\&\left.+2\fbs\int_{\Theta^K}(f(\bm y)-\fbs)\mathbf{p}_\sigma(\teb,\mathrm{d}\bm y) \right) \mub_\sigma(\mathrm{d} \teb)\\
	=&\underbrace{\int_{\Theta^K} \left( \int_{\Theta^K} (f(\bm y)-\fbs)\mathbf{p}_\sigma(\teb,\mathrm{d}\bm y)\right)^2\mub_\sigma(\mathrm{d}\teb)}_\text{$I$}+(\fbs)^2  \\ &+ 2\fbs \underbrace{\int_{\Theta^K}\int_{\Theta^K}(f(\bm y)-\fbs)\mathbf{p_\sigma(\teb,\mathrm{d}\bm y )\mub_\sigma(\mathrm{d}\teb)}}_\text{$=0$ by stationarity} \label{eq:p2}
	\end{align}
	\rjp{Thus, multiplying and dividing $I$ by $$\left(\int_{\Theta^K} \left(f(\teb)-\fbs\right)^2\mub_\sigma(\mathrm{d}\teb)\right),$$ we obtain from the definition of $\lno \mathbf{P}_\sigma \rno_{L^0_2}^2$ that:}
	
	\begin{align}
	\eqref{eq:p2}=&\left(\frac{\int_{\Theta^K} \left( \int_{\Theta^K} (f(\bm y)-\fbs)\mathbf{p}_\sigma(\teb,\mathrm{d}\bm y)\right)^2\mub_\sigma(\mathrm{d}\teb)}{ \int_{\Theta^K} \left(f(\teb)-\fbs\right)^2\mub_\sigma(\mathrm{d}\teb)}\right)\\ &\times\left(\int_{\Theta^K} \left(f(\teb)-\fbs\right)^2\mub_\sigma(\mathrm{d}\teb)\right)+(\fbs)^2 \\ 
	\leq&\lno \mathbf{P}_\sigma \rno_{L^0_2}^2\left(\int_{\Theta^K} \left(f(\teb)-\fbs\right)^2\mub_\sigma(\mathrm{d}\teb)\right)+(\fbs)^2\\
	=&\lno \mathbf{P}_\sigma \rno_{L^0_2}^2\left(\int_{\Theta^K}f(\teb)^2\mub_\sigma(\mathrm{d}\teb)\right)\\&+\left(1-\lno \mathbf{P}_\sigma \rno_{L^0_2}^2\right)(\fbs)^2 \\
	=&\left(\int_{\Theta^K}f(\teb)^2\mub_\sigma(\mathrm{d}\teb)\right)-\underbrace{\left(1-\lno \mathbf{P}_\sigma \rno_{L^0_2}^2\right)}_\text{$:=\gamma$, with $\gamma\in(0,1)$}\\ & \times\left(\int_{\Theta^K} \left(f(\teb)-\fbs\right)^2\mub_\sigma(\mathrm{d}\teb)\right) \label{Eq:boun_on_sg}.
	\end{align}
	Replacing Equation (\ref{Eq:boun_on_sg}) into Equation (\ref{Eq:ine_sg_weighted}), we get 
	\begin{align}
	&\lno \mathbf{P}^\mathrm{(W)}\rno^2_{L_2^0(\Theta^K,\mubis)\mapsto L_2^0(\Theta^K,\mubis)}
	\\&\leq \underset{f\in \mathcal{L}}{\sup}\left(\int_{\Theta^K}f(\teb)^2\mubis(\mathrm{d}\teb)\right)-\frac{\gamma}{|S_K|}\sum_{\sigma\in S_K}\V_{\mub_\sigma}[f]\\
	&\leq 1-\gamma \left(\frac{\Lambda_m}{2-\Lambda_m}\right)<1 \quad \text{(by Proposition \ref{proposition:bound_var})} \label{eq:thm_last}.\end{align}

	
	%
	
	Thus, $\mathbf{\Pw}$ has an $L_2(\Theta^K,\mubis)$ spectral gap. Once again, $L_r(\Theta^K,\mubis)$-geometric ergodicity (with $r\in[1,\infty])$ follows from Proposition \ref{proposition:l2sg_rge} and the fact that $\mathbf{P}^\mathrm{(W)}$ is $\mubis$-reversible by Proposition \ref{proposition:rev}.

\end{proof}
	
\end{thm}

	\subsubsection{Discussion and comparison to similar theoretical result} Theorems \ref{proposition:uwIS_ergodicity} and \ref{proposition:wIS_ergodicity} state the existence of an $L_2$-spectral gap, hence $L_r$-geometric ergodicity for both the UGPT and the WGPT algorithm. Their proof provides also a quantification of the $L_2$-spectral gap in terms of the $L_2$-spectral gap of each individual Markov operator $P_k$. Such a bound is, however, not satisfactory as it does not use any information on the temperature and it just states that the $L_2$-spectral gap of the UWPT and WGPT chain is not worse that the smallest $L_2$-spectral gap among the individual chains (without swapping). This result is not sharp, as it can be evidenced in the numerical section, where a substantial improvement in convergence is achieved by our methods. 
		
	Convergence results  for  the \textit{standard} parallel tempering algorithm have been obtained in the works \cite{miasojedow2013adaptive} and \cite{woodard2009conditions}. In particular, the work \cite{miasojedow2013adaptive} has proved geometric ergodicity for the  pairwise parallel tempering algorithm using the standard drift condition construction of \cite{meyn2012markov}. It is unclear from that work which convergence rate is obtained for the whole algorithm. In comparison, our results are given in terms of spectral gaps. On the other hand, the work \cite{woodard2009conditions} presents conditions for rapid mixing of a particular type of parallel tempering algorithm, where the transition kernel is to be understood as a convex combination of such kernels, as opposed to our case, where it is to be understood as  a tensorization. Their obtained results provide, for their setting, a better convergence rate that the one we obtained for the UGPT. We believe that their result can be extended to the UGPT algorithm, and this will be the focus of future work. On the other hand, the use of the ideas in \cite{woodard2009conditions} for the WGPT algorithm seems more problematic.

\section{Numerical experiments}\label{S:Num_Exp}
We now present four academic examples to illustrate the efficiency of both GPT algorithms discussed herein and compare them to the more traditional random walk Metropolis and standard PT algorithms. Notice that we compare \rjp{the different algorithms in their simplest version that uses random walk Metropolis as a base transition kernel. The only exception is in Section \ref{ss:high_dim}, which presents a high-dimensional BIP for which the preconditioned Crank-Nicolson \cite{cotter2013mcmc} is used as the base method in all algorithms instead of RWM. More advanced samplers, such as Adaptive metropolis \cite{haario2006dram,haario2001adaptive}, or other transition kernels, could be used as well to replace RWM or pCN}. Experiments \ref{ss:exp_circle}, \ref{ss:exp_bip} and \ref{ss:exp_cauchy_ac} were run in a Dell (R) Precision (TM) T3620 workstation with Intel(R) Core(TM) i7-7700 CPU  with 32 GB of  RAM. Numerical simulations in Section \ref{ss:exp_circle} and \ref{ss:exp_cauchy_ac} were run on a single thread, while the numerical simulations in Section \ref{ss:exp_bip}  were run on an \textit{embarrassingly parallel}  fashion over 8 threads using the Message Passing Interface (MPI) and the Python package MPI4py \cite{dalcin2005mpi}. Lastly, experiments \ref{ss:exp_bip_ac} and \ref{ss:high_dim} were run on the Fidis cluster of the EPFL. The scripts used to generate the results presented in this section were written in Python 3.6, and can be found in \href{https://doi.org/10.5281/zenodo.3700048}{\texttt{DOI: 10.5281/zenodo.4736623}}

\subsection{Implementation remarks}

In most Bayesian inverse problems, particularly those dealing with large-scale computational models, the computational cost is dominated by the evaluation of the  forward operator, which can be, for example, the numerical approximation of a possibly non-linear partial differential equation. In the case where all possible permutations are considered (i.e., $S_K=\mathscr{S}_K$), there are $K!$ possible permutations of the states, the computation of the swapping ratio in the GPT algorithms can become prohibitively expensive if one is to evaluate $K!$ forward models, even for moderate values of $K$. This problem can be circumvented by storing the values $\pi(\te^{(n)}_k)$, $k=1,\dots,K,$ $n=1,\dots N$, since the swapping ratio for GPT consists of permutations of these values, divided by the temperature parameters.  Thus, ``only'' $K$ forward model evaluations need to be computed at each step and the swapping ratio can be computed at negligible cost for moderate values of $K$.

\rjp{

There is, however, a clear trade-off between the choice of $K$ (which has a direct impact on the efficiency of the method), and the computational cost associated to (G)PT. Intuitively, a large $K$ would provide a better mixing, however, it requires a larger number of forward model evaluations, which tends to be costly. We remark that such a trade-off between efficiency and number of function evaluations is also present in some advanced MCMC methods, such as Hamiltonian Monte Carlo, where one needs to choose a number of time steps for the time integration (see, e.g., \cite{beskos2017geometric}). Furthermore, there is an additional constraint when choosing $S_K=\mathscr{S}_K$, and it is the \textit{permutation cost} associated to computing $r^\mathrm{{(UW)}}(\teb,\sigma)$ and $w_\sigma(\teb)$. In particular, the computation of either of those quantities has a complexity of $K!$ thus, this cost will  eventually surpass the cost of evaluating the forward model $K$ times. This is illustrated in Figure \ref{fig:costcomparisson}, where we plot the cost per sample of two different posteriors  vs $K$. These posteriors are taken from the numerical examples in Sections \ref{ss:exp_cauchy_ac} and \ref{ss:high_dim}. The posterior in Section \ref{ss:exp_cauchy_ac} is rather inexpensive to evaluate, since one can compute the forward map $\eff$ analytically (the difficulty associated to sampling from that posterior comes from its high multi-modality). On the contrary, evaluating the posterior in Section \ref{ss:high_dim} requires numerically approximating the solution to a time-dependent, second-order PDE, and as such, evaluating such a posterior is costly.  As we can see for $K\leq6$, the computational cost in both cases is dominated by the forward model evaluation. Notice that for $K\leq9$, the cost per sample from posterior \eqref{Eq:ac_post2} is dominated by the evaluation of the forward model. 

\begin{figure}
	\centering
\begin{tikzpicture}

\definecolor{color0}{rgb}{0.886274509803922,0.290196078431373,0.2}
\definecolor{color1}{rgb}{0.203921568627451,0.541176470588235,0.741176470588235}
\definecolor{color2}{rgb}{0.596078431372549,0.556862745098039,0.835294117647059}
\definecolor{color3}{rgb}{0,0,0}
	\begin{semilogyaxis}[
legend cell align={left},
legend style={
	fill opacity=0.8,
	draw opacity=1,
	text opacity=1,
	at={(0.97,0.03)},
	anchor=south east,
	draw=white!80!black
},
tick align=outside,
tick pos=left,
xmajorgrids,
xminorgrids,
yminorticks,
xmin=0.5, xmax=11.5,
xtick style={color=white!33.3333333333333!black},
ymajorgrids,
yminorgrids,
ymin=9.20778419465741e-07, ymax=100.1,
ytick style={color=white!33.3333333333333!black},
ylabel={Time},
xlabel={$K$}]
\addplot [thick, color0]
table {%
	1 1.47742557525635
	2 2.9548511505127
	3 4.43227672576904
	4 5.90970230102539
	5 7.38712787628174
	6 8.86455345153809
	7 10.3419790267944
	8 11.8194046020508
	9 13.2968301773071
	10 14.7742557525635
	11 16.2516813278198
};
\addlegendentry{Cost model \S \ref{ss:high_dim}}
\addplot [thick, color1]
table {%
	1 0.0011
	2 0.0022
	3 0.0033
	4 0.0044
	5 0.0055
	6 0.0066
	7 0.0077
	8 0.0088
	9 0.0099
	10 0.011
	11 0.0121
};
\addlegendentry{Cost model \S  \ref{ss:exp_cauchy_ac}}
\addplot [thick, color2]
table {%
	1 6.36577606201172e-05
	2 3.76701354980469e-05
	3 3.52859497070312e-05
	4 4.41074371337891e-05
	5 9.60826873779297e-05
	6 0.000436305999755859
	7 0.00319290161132812
	8 0.0291540622711182
	9 0.271856546401978
	10 2.92558073997498
	11 33.3677287101746
};
\addlegendentry{Permutation cost }

\addplot [thick, color3]
table {%
	1 0.00116365776062012
	2 0.00223767013549805
	3 0.00333528594970703
	4 0.00444410743713379
	5 0.00559608268737793
	6 0.00703630599975586
	7 0.0108929016113281
	8 0.0379540622711182
	9 0.281756546401978
	10 2.93658073997498
	11 33.3798287101746
};
\addlegendentry{Total cost \S \ref{ss:high_dim}}

\addplot [thick, white!46.6666666666667!black]
table {%
	1 1.47748923301697
	2 2.95488882064819
	3 4.43231201171875
	4 5.90974640846252
	5 7.38722395896912
	6 8.86498975753784
	7 10.3451719284058
	8 11.8485586643219
	9 13.5686867237091
	10 17.6998364925385
	11 49.6194100379944
};
\addlegendentry{Total cost \S \ref{ss:exp_cauchy_ac}}

\end{semilogyaxis}%

\end{tikzpicture}
	\caption{Cost per sample vs $K$ for $S_K=\mathscr{S}_K$ for the forward model in Section \ref{ss:exp_cauchy_ac} and the forward model in \ref{ss:high_dim}.}
	\label{fig:costcomparisson}
\end{figure}
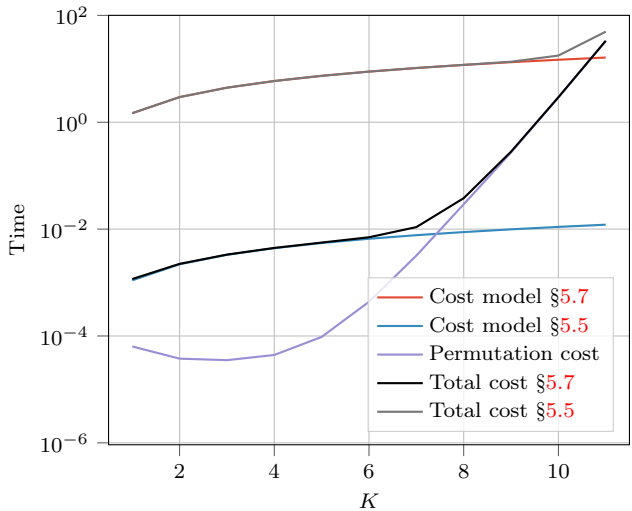
}
Thus, for high values of $K$, it is advisable to only consider the union of properly chosen  semi-groups  $A,B$ of $\mathscr{S}_K,$ with $ A\cap B\neq \emptyset,$ such that $A,B$ generates $\mathscr{S}_K$ (i.e., if the smallest semi-groups that contains $A$ and $B$ is $\mathscr{S}_K$ itself), and $|A\cup B|<|\mathscr{S}_K|=K!,$ which is referred to as partial Infinite Swapping in the continuous case \cite{dupuis2012infinite}. One particular way of choosing $A$ and $B$ is to consider, for example, $A$ to be the set of permutations that only permute the indices associated with relatively low temperatures while leaving the other indices unchanged, and $B$ as the set of permutations for the indices of relatively high temperatures, while leaving the other indices unchanged. Intuitively, swaps between temperatures that are, in a sense, ``close'' to each other tend to be chosen with a higher probability.   We refer the reader to \cite[Section 6.2]{dupuis2012infinite} for a further discussion on this approach in the continuous-time setting. One additional idea would be to consider swapping schemes that, for example, only permute states  between $\mu_i$ and $\mu_{i+1}, \mu_{i+2}, \dots, \mu_{i+\ell}$ for some user-defined $\ell \geq 1$ and any given $i=1,2,\dots, K-1$. The intuition behind this choice also being that swaps between posteriors that are at close temperatures are more likely to occur than swaps between posteriors with a high temperature difference. We intend to explore this further in depth in future work.

\rjp{We reiterate that the total number of temperatures $K$ depends heavily on the problem and the computational budget available \cite{doll2012rare,van2008parameter,Yu11744} For the experiments considered in the work we chose $K=4$ or $K=5$, which provide an acceptable  compromise between acceleration and cost.}

\rjp{

\subsection{Experimental setup}	
	
We now present an experimental setup common to all the numerical examples presented in the following subsections. In particular, all the experiments presented in this work utilize a \textit{base} method given by either RWM  (for experiments \ref{ss:exp_circle} through \ref{ss:exp_bip_ac}) or pCN (used in experiment \ref{ss:high_dim}) for the Markov transition kernels $p$. Furthermore, we take $S_K=\mathscr{S}_K$ for all experiments, where  $K=5$ for experiment \ref{ss:exp_cauchy_ac} and $K=4$ for the other 4 experiments. In addition, we follow the \textit{rule of thumb} of \cite{earl2005parallel} for the choice of temperatures, setting, for each experiment, $T_{k}=a^{k-1},$ $k=1,\dots,K$, for some positive  constant $a>1$. The particular choice of $a$ is problem-dependent and it is generally chosen so that $\mu_K$ becomes sufficiently simple to explore. For each experiment we implement 5 MCMC algorithms to sample from a given posterior $\mu=\mu_1$, namely, the base (untempered) method (either RWM or pCN), and such a method combined with the standard PT algorithm (PT) with $N_s=1,$ the PSDPT algorithm of \cite{lkacki2016state}, and both versions of  GPT. For our setting, the tempered algorithms have a cost (in terms of number of likelihood evaluations) that is $K$ times larger than the base method. Thus, to obtain a fair comparison across all algorithms, we run the chain for the base method  $K$ times longer.   Lastly, given some problem-dependent quantity of interest $\mathcal{Q}$,  we assess the efficiency of our proposed algorithms to compute the posterior expectation of $\mathcal{Q}$ by  comparing the mean square error (experiments \ref{ss:exp_circle}-\ref{ss:exp_cauchy_ac}), for which the exact value of $\E_{\mu}[\mathcal{Q}]$ is known, or the variance (experiments \ref{ss:exp_bip_ac}-\ref{ss:high_dim}) of the ergodic estimator $\hat{\mathcal{Q}}$ obtained over $N_\text{runs}$ independent runs of each algorithm. 
}

\subsection{Density concentrated over a quarter circle-shaped manifold} \label{ss:exp_circle}

\begin{figure}[tbh]
	\centering
	\includegraphics[width=1.\linewidth]{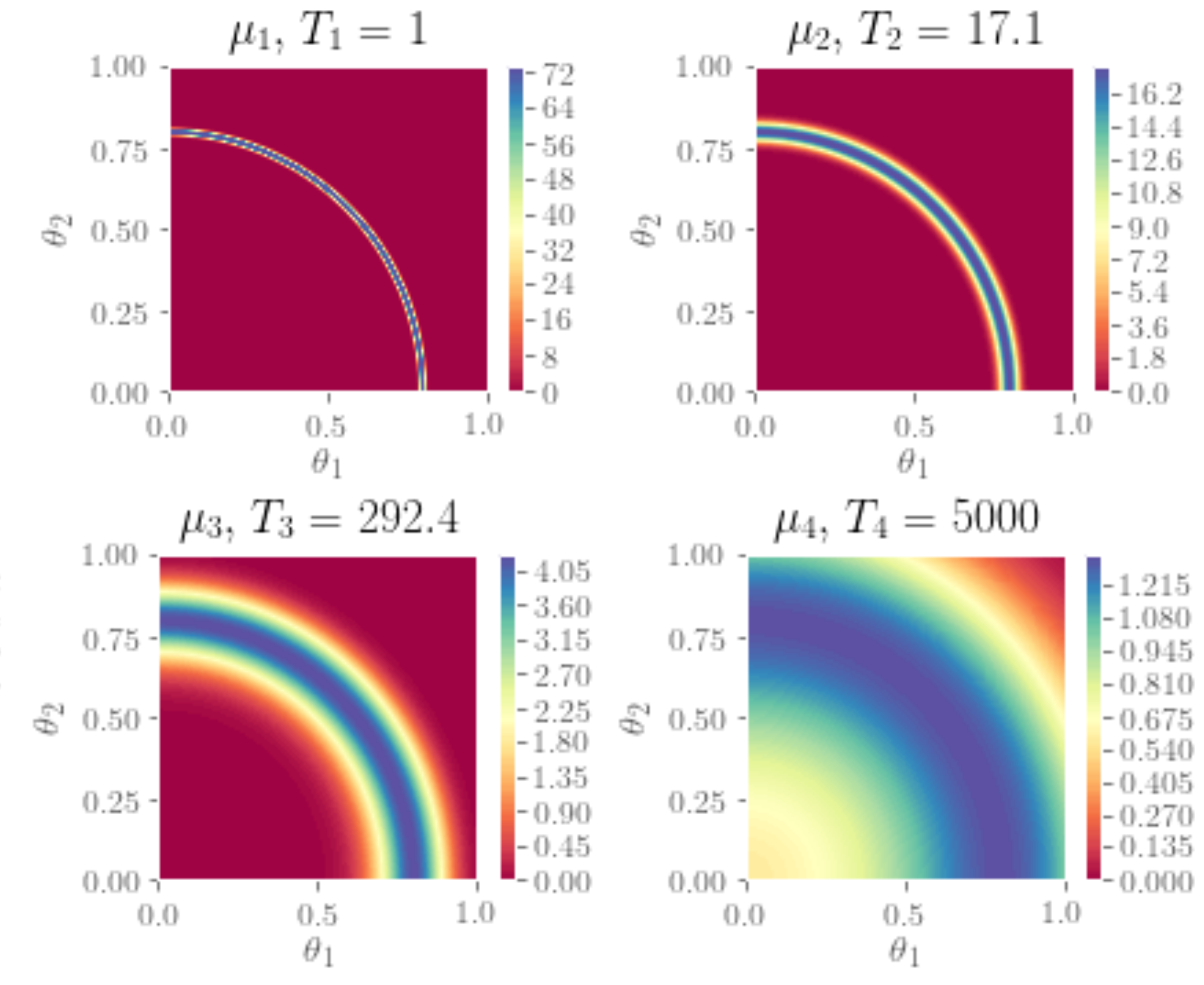}
	\caption{Tempered densities (with $T_1=1,\ T_2=17.1,\ T_3=292.4,\ T_4=5000$) for the density concentrated around a quarter circle-shaped manifold example. As we can see, the density becomes less concentrated as the temperature increases, which allows us to use RWM proposals with larger step sizes.}
	\label{fig:truecirlcedensity4}
\end{figure}

Let $\mu$ be a probability measure that has density $\pi$  with respect to the uniform Lebesgue measure on the unit square $\mu_\mathrm{pr}=\mathcal{U}([0,1]^2)$ given by \begin{align}
\pi(\te)=\frac{1}{Z}\exp\left(-10000(\te_1^2+\te_2^2-0.8^2)^2	  \right)\mathbf 1_{[0,1]^2},\label{Eq:density_circle} \
\end{align}
where $\te=(\te_1,\te_2),$ $Z$ is the normalization constant, and $\mathbf1_{[0,1]^2}$ is the indicator function over the unit square.  We remark that this example is not of particular interest \textit{per se}; however, it can be used to illustrate some of the advantages of the algorithms discussed herein. The difficulty of sampling from such a distribution comes from the fact that its density is concentrated over a quarter circle-shaped manifold, as  can be seen on  the left-most plot in Figure \ref{fig:truecirlcedensity4}. This in turn will imply that a single level RWM chain would need to take very small steps in order to properly explore such density.

 We aim at estimating  $\hat{\mathcal{Q}}_i=\E_{\mu}[\te_i]\approx \hat\te_i$, for $i=1,2$.  For the tempered algorithms (PT, PSDPT, UGPT,  and WGPT), we consider $K=4$ temperatures and choose $T_4=5000$, so that the tempered density $\pi_4$ becomes sufficiently simple to explore the target distribution. This gives $T_1=1, T_2=17.1, T_3=292.4, T_4=5000$.
 We compare the quality of our algorithms by examining the variance of the estimators $\hat\te_i,$ $i=1,2$ computed over $N_\text{runs}=100$ independent MCMC runs of each algorithm. For the tempered algorithms, each estimator is obtained by running the inversion experiment for $N=25,000$ samples per run, discarding the first 20\%  of the samples (5000) as a burn-in. Accordingly, we run the single-chain random walk Metropolis algorithm for $N_\mathrm{RWM}=K N=100, 000$ iterations, per run, and discard the first 20\% of the  samples obtained with the RWM algorithm (20,000) as a burn-in.

The untempered RWM algorithm uses Gaussian proposals with covariance matrix $\Sigma_\mathrm{RWM}=\rho_1^2 I_{2\times 2},$ where $I_{2\times2}$ is the identity matrix in $\R^{2\times 2}$, \rjp{and $\rho^2_1=0.022$ is chosen in order to obtain an acceptance rate of around 0.23}. For the tempered algorithms (i.e., PT, PSDPT, and both versions of GPT),  we use $K=4$ RWM kernels $p_k$, $k=1,2,3,4$, with proposal density $q_{\mathrm{prop},k}(\te^{(n)}_k,\cdot)=\mathcal{N}(\te_k^{(n)},\rho_k^2I_{2\times 2})$, where $\rho_k$ is shown in Table \ref{tab:sigma_banana}. This choice of $\rho_k$  gives an acceptance rate for each chain of around 0.23.  Notice that $\rho_1$ corresponds to the ``step-size'' of the single-temperature RWM algorithm.

\begin{table}[tbh]
	\centering
	\begin{tabular}{lrrrr}
\hline
		&        $k=1$ &  $k=2$ &  $k=3$ &        $k=4$  \\
\hline
		$\rho_k$     &  0.022 &  0.090 &  0.310 &  0.650 \\
\hline
	\end{tabular}
	\caption{Step size of the RWM proposal distribution for the manifold experiment. \label{tab:sigma_banana} }

\end{table}

\begin{figure}[tbh]
	\centering
		\includegraphics[width=1\linewidth]{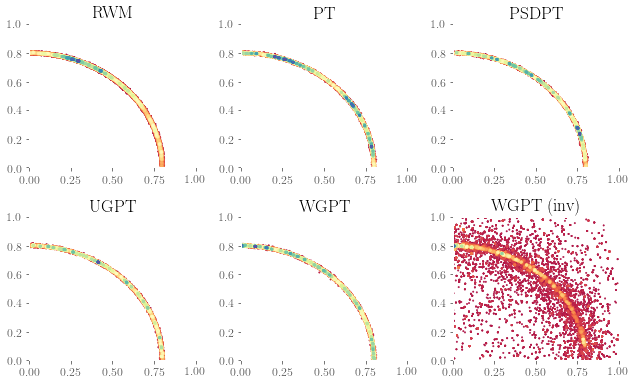}
	\caption{ Scatter-plots of the samples from $\mu$ obtained with each algorithm on a single run. Top, from left to right: random walk Metropolis, PT and PSDPT. Bottom, from left to right: UGPT, WGPT (after re-weighting the samples), and WGPT, before re-weighting the samples.}
	\label{fig:2dhistogram_all}
\end{figure}

  Experimental results for the ergodic run  are shown in Table \ref{tab:res_banana}. We can see how both GPT algorithms provide a gain over  RWM, PT and PSDPT algorithms, with the WGPT algorithm providing the largest gain. Scatter plots of the samples obtained with each method are presented in Figure \ref{fig:2dhistogram_all}. Here, the subplot titled ``WGPT'' (\rjp{bottom row, middle}) corresponds to weighted samples from $\mub_{\textrm{W}}$, with weight $\hat{w}$ as in \eqref{Eq:weighted_samples}, while the one titled  ``WGPT (inv)'' \rjp{(bottom row, right)} corresponds to samples from  $\mubis$ without any post-processing. Notice how the samples from the latter concentrates over a \textit{thicker} manifold, which in turn makes the target density easier to explore when using state-dependent Markov transition kernels.
\begin{table}[tbh]
		\centering
			\resizebox{\columnwidth}{!}{%
			\begin{tabular}{lcccccc}
				\hline
				& \multicolumn{2}{c}{Mean} & \multicolumn{2}{c}{MSE} &  \multicolumn{2}{c}{MSE$_\mathrm{RWM}/$MSE}\\
				& $\hat{\te}_1$     & $\hat{\te}_2$          & $\hat{\te}_1$     & $\hat{\te}_2$   & $\hat{\te}_1$     & $\hat{\te}_2$ \\ \hline
				\multicolumn{1}{c}{} &             &            &               &              &                          \\
RWM 	&0.50996&	0.50657&			0.00253&	0.00261&	1.00&	1.00\\
PT		&0.50978&	0.51241&			0.00024&	0.00021&	10.7&	11.0\\
PSDPT	&0.50900&	0.50956&			0.00027&	0.00026&	9.53&	10.2\\
UGPT	&0.50986&	0.50987&			0.00016&	0.00016&	16.1&	16.4\\
WGPT	&0.51062&	0.50838&			0.00015&	0.00014&	16.9&	18.4  \\\hline
		\end{tabular}}
		\caption{Results for the density concentrated around a circle-shaped manifold experiment. As we can see, both GPT algorithms provide an improvement over PT, PSDPT and RWM.  The computational cost is comparable across all algorithms. \label{tab:res_banana}}
			\label{tab:erg_banana}

\end{table}

%

\subsection{Multiple source elliptic BIP}\label{ss:exp_bip}
We now consider a slightly more challenging problem, for which we try to recover the probability distribution of  the location of a source term in a Poisson equation (Eq. (\ref{Eq:elliptic})), based on some noisy measured data. Let $(\Theta,\mathcal{B}(\Theta),\mu_\mathrm{pr})$ be the measure space, set
$\Theta=\bar D:=[0,1]^2$, with Lebesgue (uniform) measure $\mu_\mathrm{pr},$  and consider the following Poisson's equation with homogeneous boundary conditions:
 \begin{align}
\begin{cases}
\Delta u(x,\te)=f( x,\te),&  x\in D,\  \te\in \Theta,\\
u(x,\te)=0,& x\in \partial D.
\end{cases}\label{Eq:elliptic}
\end{align}
 Such equation can model, for example, the electrostatic potential $u:=u( x,\te)$  generated by a charge density $f( x,\te)$ depending on an \textit{uncertain} location parameter $\te\in \Theta$.  Data $y$ is recorded on an array of $64\times64$ equally-spaced points in $D$ by solving \eqref{Eq:elliptic}  with a forcing term  given by
 \begin{align} \label{Eq:force_elliptic}
 f(x)= \sum_{i=1}^{4}e^{-1000[(x_1-s ^{(i)}_1)^2
 	+(x_2-s^{(i)}_2)^2]},
 \end{align}
 where the true source locations $s^{(i)}, \ i=1,2,3,4,$ are given by $s^{(1)}=(0.2,0.2),\ s^{(2)}=(0.2,0.8),\ s^{(3)}=(0.8,0.2),$ and $s^{(4)}=(0.8,0.8)$. Such data is assumed to be polluted by an additive Gaussian noise with distribution $\mathcal{N}(0,\eta^2 I_{64\times64})$, with $\eta=3.2\times10^{-6}$, (which corresponds to a $1$\% noise) and where $I_{64\times64}$ is the 64-dimensional identity matrix. Thus, we set $(Y,\lno \cdot\rno_Y)=(\R^{64\times64},\lno \cdot\rno)$, with  $\lno A \rno=(64\eta)^{-2}\lno A \rno^2_F $, for some arbitrary matrix $A\in \R^{64\times 64},$ where $\lno \cdot \rno_F$ is the Frobenius norm.
We assume a misspecified model where we only consider a single source in Eq. \eqref{Eq:force_elliptic}. That, is, we construct our forward operator $\eff: \Theta\mapsto Y$ by solving \eqref{Eq:elliptic} with a source term given by   \begin{align} \label{Eq:forcing_term}
 f(x,\te)= e^{-1000[(x_1-\te_1)^2 +(x_2-\te_2)^2]}.
 \end{align}
In this particular setting, this leads to a posterior distribution with four modes since the prior density is uniform in the domain and the likelihood has a local maximum whenever $(\te_1,\te_2)=(s_1^{(i)},s_2^{(i)}),\ i=1,2,3,4$.  The Bayesian inverse problem at hand can be understood as sampling from the posterior measure $\mu$, which has a density with respect to the prior $\mu_\mathrm{pr}=\mathcal{U}(\bar D)$ given by \begin{align}\label{Eq:multi_modal_source}
\pi(\te)=\frac{1}{Z}\exp\left(-\frac{1}{2}\lno y-\mathcal{F}(\te)\rno^2_\Sigma\right),
\end{align}
for some (intractable) normalization constant $Z$ as in \eqref{Eq:BT}. We remark that the  solution to (\ref{Eq:elliptic}) with a forcing term of the form of \eqref{Eq:forcing_term} is approximated using a second-order accurate finite difference approximation with grid-size $h=1/64$ on each spatial component.

\begin{figure}
	\centering
	\includegraphics[width=1\linewidth]{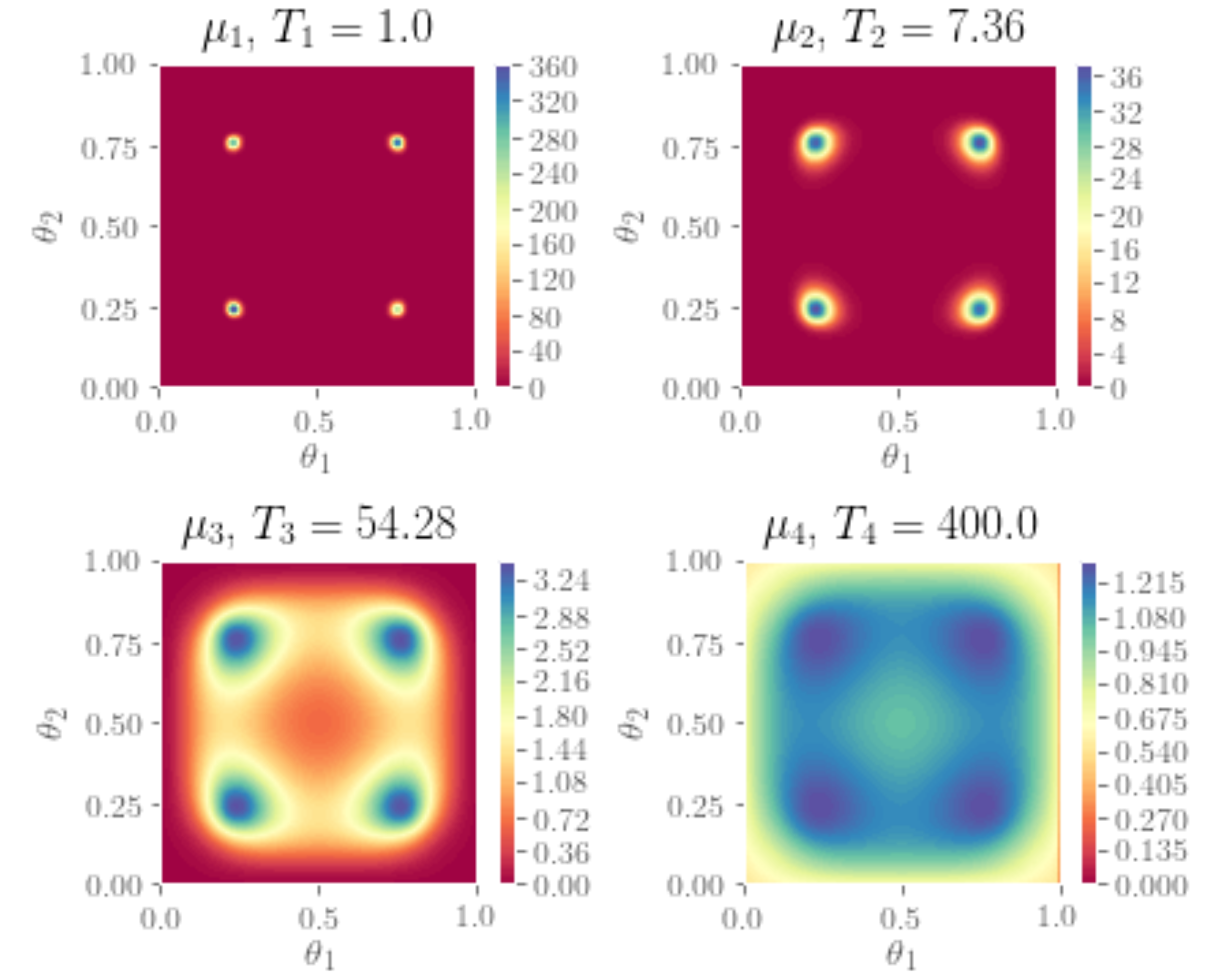}
	\caption{True tempered densities for the elliptic BIP example. Notice that the density is not symmetric, due to the additional random noise.}
	\label{fig:trueellipticdensity4}
\end{figure}

The difficulty in sampling from the current BIP arises from the fact that the resulting posterior $\mu$ is multi-modal and the number of modes is not known apriori (see Figure \ref{fig:trueellipticdensity4}).

We follow a similar experimental setup  to the previous example, and aim at estimating  $\hat{\mathcal{Q}}_i=\E_{\mu}[\te_i]\approx \hat\te_i$, for $i=1,2$. 
  We use $K=4$ temperatures and $N_\text{runs}=100$. For the PT, PSDPT and GPT algorithms, four different temperatures are used, with $T_1=1,\ \ T_2=7.36,\ \ T_3=54.28,$ and  $T_4=400$.  For each run, we obtain $N=25,000$ samples with the PT, PSDPT, and both GPT algorithms, and $N=100,000$ samples with RWM, discarding the first 20\% of the samples in both cases (5000, 20000, respectively) as a burn-in.  On each of the tempered chains, we use RWM proposals, with step-sizes shown in table \ref{tab:size_elliptic}. This choice of step size provides an acceptance rate of about 0.24 across all tempered chains and all tempered algorithms. For the single-temperature RWM run, we choose a larger step size ($\rho_\mathrm{RWM}=0.16$) so that the RWM algorithm is able to explore the whole distribution. Such a choice, however, provides a smaller acceptance rate of about 0.01 for the single-chain RWM.

    Experimental results are shown in Table \ref{tab:res_elliptic}. Once again, we can see how both GPT algorithms provide a gain over RWM and both variations of the PT algorithm, with the WGPT algorithm providing a larger gain. Scatter-plots of the obtained samples are shown in Figure \ref{fig:trueellipticdensity4}.

  \begin{table}[tbh]
  	\centering
  	\begin{tabular}{lrrrr}
\hline
  		&        $k=1$ &  $k=2$ &  $k=3$ &        $k=4$  \\
\hline
  		$\rho_{k,\mathrm{Tempered}}$     &  0.030 &  0.100 &  0.400 &  0.600 \\
  	  	$\rho_{k,\mathrm{RWM}}$     &  0.160 &  - &  - &  - \\
\hline
  	\end{tabular}
  	\caption{Step size of the RWM proposal distribution for the elliptic BIP experiment.  	\label{tab:size_elliptic} }

  \end{table}
\begin{table}[tbh]
\centering
			\resizebox{\columnwidth}{!}{%
	\begin{tabular}{lcccccc}
		\hline
		& \multicolumn{2}{c}{Mean} & \multicolumn{2}{c}{MSE} &  \multicolumn{2}{c}{MSE$_\mathrm{RWM}/$MSE}\\
		& $\te_1$     & $\te_2$          & $\te_1$     & $\te_2$   & $\te_1$     & $\te_2$ \\ \hline
		\multicolumn{1}{c}{} &             &            &               &              &                          \\
		RWM                  & 0.48509     & 0.51867       & 0.00986    & 0.01270    &1.00&  1.00   \\
		PT              	    & 0.48731     & 0.50758        & 0.00042    & 0.00036    &23.0& 29.2    \\
		PSDPT		         & 0.48401     & 0.50542       & 0.00079    & 0.00099   &12.4& 10.7    \\
		UGPT                & 0.48624     & 0.50620        & 0.00038    & 0.00027    &25.9& 38.2   \\
		WGPT               & 0.48617     & 0.50554        & 0.00025    & 0.00023    &38.6& 44.9     \\\hline
\end{tabular}}

	\caption{Results for the elliptic BIP problem.  The computational cost is comparable across all algorithms, given that the cost of each iteration is dominated by the cost of solving the underlying PDE. \label{tab:res_elliptic} }

\end{table}

%

\begin{figure}
	\centering
	\includegraphics[width=1\linewidth]{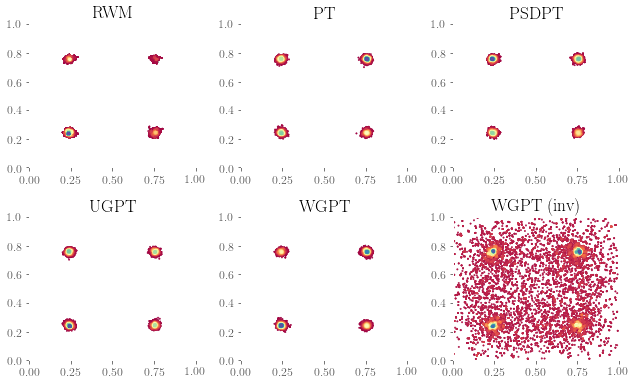}
	\caption{Scatterplots of the samples from $\mu$ obtained with different algorithms on a single run. Top, from left to right: random walk Metropolis, PT and PSDPT. Bottom, from left to right: UGPT, WGPT (after re-weighting the samples), and WGPT, before re-weighting the samples. As we can see,  WGPT (before re-weighting) is able to "connect" the parameter space.}
	\label{fig:allscatterplotselliptic}
\end{figure}

\subsection{1D wave source inversion}\label{ss:exp_cauchy_ac}
\begin{figure}
	\centering
	\input{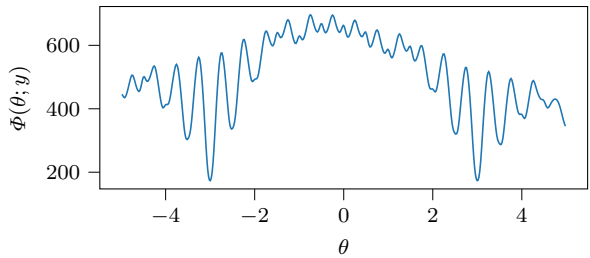}
	\caption{Multi-modal potential for the Cauchy problem. Notice the  minima around $\te=-3$ and $\te=3$.  }
	\label{fig:nll}
\end{figure}

We consider a small variation of example 5.1 in \cite{motamed2019wasserstein}.  Let $(\Theta,\mathcal{B}(\Theta),\mu_\mathrm{pr})$ be a measure space, with $\Theta=[-5,5]$ and uniform (Lebesgue) measure $\mu_\mathrm{pr}$, and let $I=(0,T]$ be a time interval.   Consider the following Cauchy problem for the 1D wave equation:
\begin{align}
\begin{cases}
u_{tt}(x,t,\te)-u_{xx}(x,t,\te)=0,& (x,t,\theta) \in \R\times I\times \Theta,\\
u(x,0,\te)=h(x,\te), &(x,t,\theta) \in \R\times\{0\}\times\Theta,\\
u_t(x,0,\te)=0,& (x,t,\theta) \in \R\times\{0\}\times\Theta. \\
\end{cases}\label{Eq:cauchy}
\end{align}
Here, $h(x,\te)$ acts as a source term generating a initial wave pulse. Notice that Equation \eqref{Eq:cauchy} can be easily solved using d'Alembert's formula, namely
\begin{align}
u(x,t,\te)=\frac{1}{2}\left(h(x-t,\te)+h(x+t,\te)\right).
\end{align}
Synthetic data $y$ is generated by solving Equation \eqref{Eq:cauchy} with initial data\begin{align}\label{eq:initial_data}
 h(x,\te_1,\te_2)&=\frac{1}{2}\left(e^{-100(x-\te_1-0.5)^2}
 +e^{-100(x-\te_1)^2}\right.\\
 &\left.+e^{-100(x-\te_1+0.5)^2}+e^{-100(x-\te_2-0.5)^2}\right. \\ &\left.+e^{-100(x-\te_2)^2}+e^{-100(x-\te_2+0.5)^2}\right),
\end{align}
with $\te_1=-3,\te_2=3$ and observed at $N_R=11$ equally-spaced receiver locations between $R_1=-5$ and $R_2=5$ on $N_T=1000$ time instants between $t=0$ and $T=5$. The signal recorded by each receiver is assumed to be polluted by additive Gaussian noise $\mathcal{N}(0,\eta^2 I_{1000\times1000})$, with $\eta=0.01$, which corresponds to roughly 1\% noise.
 We set $(Y,\lno  \rno_{Y})=(\R^{11\times 1000}, \lno \cdot\rno_{\Sigma}), $ with $$\lno A \rno_{\Sigma}^2=(\sqrt{N_R}\eta)^{-2}\sum_{i=1}^{N_R}\sum_{j=1}^{N_T}A_{i,j}^2,$$ $A\in\R^{11\times1000}$. Once again, we assume a misspecified model where we  construct our forward operator $\eff: \Theta\mapsto Y$ by solving \eqref{Eq:cauchy} with a source term given by   \begin{align} \label{Eq:h}
h(x,\te)&=\left(e^{-100(x-\te-0.5)^2}\right. \\ &\left.+e^{-100(x-\te)^2}+e^{-100(x-\te+0.5)^2}\right).
\end{align}
The Bayesian inverse problem at hand can be understood as sampling from the posterior measure $\mu$, which has a density with respect to the prior $\mu_\mathrm{pr}=\mathcal{U}([-5,5])$ given by \begin{align}\label{Eq:cauchy_b}
\pi(\te)&=\frac{1}{Z}\exp\left(-\frac{1}{2}\lno y-\mathcal{F}(\te)\rno^2_\Sigma\right)\\
&=\frac{1}{Z}\exp\left(-\Phi(\te;y)\right),
\end{align}
for some (intractable) normalization constant $Z$ as in \eqref{Eq:BT}. The difficulty in solving this BIP comes from the high multi-modality of the potential $\pot$, as it can be seen in Figure \ref{fig:nll}. This, shape of $\pot$ makes the posterior difficult to explore using local proposals.

 In this case, we consider $K=5$, and set $T_1=1,\ \ T_2=5,\ \ T_3=25,$ $T_4=125$ and $T_5=625$. \rjp{Notice that from Figure \ref{fig:costcomparisson}, the computational cost per sample is dominated by the evaluation of \eqref{Eq:cauchy_b} for values of $K\leq 6$. } Once again,  we obtain $N=25,000$ samples with the PT, PSDPT, and both GPT algorithms, and $N=125,000$ samples with RWM, discarding the first 20\% of the samples in both cases (5000, 25000, respectively) as a burn-in.  On each of the tempered chains, we use RWM proposals, with step-sizes shown in table \ref{tab:size_cauchy}. This choice of step size provides an acceptance rate of about 0.4 across all tempered chains and all tempered algorithms. The choice of step-size for the RWM algorithm is done in such a way that it can "jump" modes, which are at distance of roughly 1/2. 

We consider $\mathcal{Q}=\te$ as a quantity of interest. Experimental results are shown in Table \ref{tab:res_cauchy}. Once again, we can see how both GPT algorithms provide a gain over RWM and both variations of the PT algorithm, with the WGPT algorithm providing the largest gain. Notice that, given the high muti-modality of the posterior at hand, the simple RWM algorithm is not well-suited for this type of distribution, as it can be seen from its large variance; this suggests that the RWM usually gets "stuck"  at one mode of the posterior. Notice that, intuitively, due to the symmetric nature of the potential, one would expect the true mean of $\te$ to be close to 0. This value was computed by means of numerical integration and is given by $\E_\mu[\te]=0.08211$.
\begin{table}[tbh]
	\centering
	\begin{tabular}{lrrrrr}
		\hline
		&        $k=1$ &  $k=2$ &  $k=3$ &        $k=4$  & $k=5$ \\
		\hline
		$\rho_{k,\mathrm{Tempered}}$     &  0.02 &  0.05 &  0.10 &  0.50 & 2.0\\
		$\rho_{k,\mathrm{RWM}}$     &  0.5 &  - &  - &  - &- \\
		\hline
	\end{tabular}
	\caption{Step size of the RWM proposal distribution for the Cauchy BIP experiment. 	\label{tab:size_cauchy}  }

\end{table}

\begin{table}[tbh]

	\centering
	\small{
		\begin{tabular}{lccc}
			\hline
			& \multicolumn{1}{c}{Mean} & \multicolumn{1}{c}{MSE} &  \multicolumn{1}{c}{MSE$_\mathrm{RWM}/$MSE}\\\hline
			RWM                  & -0.10120     & 9.36709   &  1.000   \\
			PT              	  & 0.05118     & 0.03681   & 254.5    \\
			PSDPT		         &  0.15840     & 0.21701   & 43.20    \\
			UGPT                &   0.08976     & 0.03032   & 308.9   \\
			WGPT               &    0.06149     & 0.02518   & 372.0     \\\hline
	\end{tabular}}

	\caption{Results for the 1D Cauchy BIP problem.  The computational cost is comparable across all algorithms.}
	\label{tab:res_cauchy}

\end{table}
\subsection{Acoustic wave source inversion}\label{ss:exp_bip_ac}
We consider a  more challenging problem, for which we try to recover the probability distribution of  the spatial location of a (point-like) source term, together with the material properties of the medium,  on an acoustic wave equation (Eq. (\ref{Eq:acoustic})), based on some noisy measured data. We begin by describing the mathematical model of such wave phenomena. Let $(\Theta,\mathcal{B}(\Theta),\mu_\mathrm{pr})$ be the measure space , with Lebesgue (uniform) measure $\mu_\mathrm{pr}$, set
$\bar D:=[0,3]\times[0,2],$ $\partial D=\bar\Gamma_{N}\cup \bar\Gamma_\mathrm{Abs},\ \mathring\Gamma_\mathrm{N}\cap \mathring \Gamma_\mathrm{Abs}=0,$ $|\Gamma_\mathrm{N}|,  |\Gamma_\mathrm{Abs}|>0,$  and define  $\Theta=D\times\Theta_\alpha\times\Theta_\beta,$ where $\Theta_\alpha=[6,14]$, $\Theta_\beta=[4500,5500]$. Here, we are considering a rectangular spatial domain $D$, with the top boundary denoted by $\Gamma_\mathrm{N}$ and the side and bottom boundaries denoted by $\Gamma_\mathrm{Abs}$. Lastly, let $\te:=(s_1,s_2,\alpha,\beta)\in\Theta$. Consider the following acoustic wave equation with absorbing boundary conditions:
\begin{align}
\begin{cases}
\alpha^2 u_{tt}-\nabla \cdot(\beta^2\nabla u)=f,&  \text{in } D\times(0,T)\times \Theta,\\
u=u_t=0,&  \text{in } D\times\{0\}\times \Theta,\\
\beta^2 \nabla u\cdot \hat n=0,&  \text{on }\Gamma_\mathrm{N}\times(0,T)\times \Theta,\\
\beta^2 \nabla u\cdot \hat n=-\alpha\beta u_t,&  \text{on } \Gamma_\mathrm{Abs}\times(0,T)\times \Theta,
\end{cases}\label{Eq:acoustic}
\end{align}
where $u=u(x,t,\te)$, and $f=f(x,t,\theta)$. Here the boundary condition on the top boundary $\Gamma_\mathrm{N}$ corresponds to a Neumann boundary condition, while the boundary condition on $\Gamma_\mathrm{Abs}$ corresponds to the so-called absorbing boundary condition, a  type of artificial boundary condition  used to minimize reflection of wave hitting the boundary.  Data $y\in Y$ is obtained by solving Equation $\eqref{Eq:acoustic}$ with a force term given by

	  \begin{align}\label{eq:force}
f(x,t,\te)&=10^{11} e^{-\frac{1}{2\cdot 0.1^2}\left[\left(x_1-s_1\right)^2+\left(x_2-s_2\right)^2\right]}\\&\times(1-2\cdot 1000\pi^2t^2)e^{-2\cdot 1000^2\pi^2t^2},
\end{align}
with a true set of parameters $\Theta\ni\te^*:=(s_1,s_2,\alpha,\beta)$ given by $s_1=1.5,s_2=1.0,\alpha=10,\beta=5000$, and observed on $N_R=3$ different receiver locations $R_1=(1.0,2.0),R_2=(1.5,2.0),R_3=(2.0,2.0)$ at $N_T=117$ equally-spaced time instants between $t=0$ and $t=0.004$. In physical terms, the parameters $s_1,s_2$ represent the source location, while the parameters $\alpha,\beta$ are related to the material properties of the medium. Notice that, on a slight abuse of notation, we have used the symbol $\pi$ to represent the number $3.14159\dots$ in equation \eqref{eq:force} and it should not be confused with the symbol for density. The data measured by each receiver is polluted by additive Gaussian noise $\mathcal{N}(0,\eta^2 I_{117\times117})$, with $\eta=0.013$, which corresponds to roughly a 2\% noise. Thus, we have that $(Y,\lno \cdot\rno_Y)=(\R^{3\times 117},\lno \cdot\rno_\Sigma), $ where $\lno A\rno_\Sigma^2:=(\sqrt{N_R}\eta)^{-2}\sum_{i=1}^{N_R}\sum_{j=0}^{N_T}A_{i,j}^2$. Thus, the forward mapping operator $\mathcal{F}:\Theta\mapsto Y$ can be understood as the numerical solution of Equation \eqref{Eq:acoustic} evaluated at 117 discrete time instants at each of the 3 receiver locations. Such a numerical approximation is obtained by the finite element method using piece-wise linear elements and the time stepping is done using  a Forward Euler scheme with  sufficiently small time-steps to respect the so-called Courant-Friedrichs-Lewy condition \cite{quarteroni2009numerical}. This numerical solution is implemented using the Python library FEniCS \cite{logg2012automated}, using 40$\times 40$ triangular elements. The Bayesian inverse problem at hand can thus be understood as sampling from the posterior measure $\mu$, which has a density with respect to the prior $\mu_\mathrm{pr}=\mathcal{U}(\Theta)$ given by \begin{align}\label{Eq:ac_post}
\pi(\te)=\frac{1}{Z}\exp\left(-\frac{1}{2}\lno y-\mathcal{F}(\te)\rno^2_\Sigma\right).
\end{align}

\begin{figure}
	\centering
	\includegraphics[width=1\linewidth]{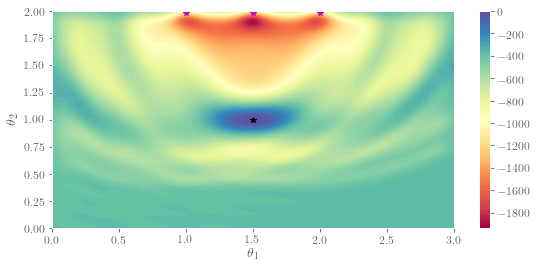}
	\caption{Plot of the log-likelihood for different values of $s_1$, $s_2$ and fixed values of $\alpha=10$ and $\beta=5000$. The magenta points represent the reciever locations $R_1,R_2,R_3$. The black point represents the true location of the source $(s_1,s_2)=(1.5,1.0)$.  }
	\label{fig:nllac}
\end{figure}

The previous BIP presents two difficulties; on the one hand,  Equation \eqref{Eq:acoustic} is, typically, expensive to solve, which in turn translates into expensive evaluations of the posterior density. On the other, the  log-likelihood has an extremely complicated structure, which in turn makes its exploration difficult. This can be seen in Figure \ref{fig:nllac}, where we  plot of the log-likelihood for different source locations $(s_1,s_2)$ and for fixed values of the material properties $\alpha=10,\beta=5000. $ More precisely, we plot  $\tilde \Phi((s_1,s_2);y):=-\frac{1}{2}\lno y-\eff(s_1,s_2,10,5000)\rno^2_\Sigma$ on a grid of $100\times 100$ equally spaced points $(s_1,s_2)$ in $D$. It can be seen that, even though the  log-likelihood has a clear peak around the true value of $(s_1,s_2)$, there are also regions of relatively high concentration of log-probability, surrounded by regions with a significantly smaller log-probability, making it a suitable problem for our setting.

Following the same set-up of previous experiments, we aim at estimating  $\hat{\mathcal{Q}}_i=\E_{\mu}[\te_i]\approx \hat\te_i$, for $i=1,2$.  Once again, we consider $K=4$ temperatures for the tempered algorithms (PT, PSDPT, UGPT,  and WGPT), and set  temperatures to  $T_1=1, T_2=7.36, T_3=54.28,$  $T_4=400$.  We compare the quality of our algorithms by examining the variance of the estimators $\hat\te_i,$ $i=1,2$ computed over $N_\text{runs}=50$ independent MCMC runs of each algorithm. For each run, we run the tempered algorithms obtaining  $N=7,000$ samples, discarding the first 20\%  of the samples (1400) as a burn-in. For the RWM algorithm, we run the inversion experiment for  $N_\mathrm{RWM}=K N=28, 000$ iterations, and discard the first 20\% of the  samples obtained (5600) as a burn-in.

Each individual chain is constructed using Gaussian RWM proposals $q_{\mathrm{prop},k}(\te^{n}_k,\cdot)=\mathcal{N}(\te^{n}_k,\mathcal{C}_k)$, $k=1,2,3,4$, with covariance $\mathcal{C}_k$ described in Table \ref{tab:size_ac}. The covariance is tuned in such a way that the acceptance rate of each chain is around 0.2.  The variance of the estimators obtained with each method is presented in Table \ref{tab:res_ac}. Once again, our GPT algorithms  outperform all other tested methods for this particular setting. In particular, our methods provide huge computational gains when compared to RWM and the PSDPT algorithm of \cite{lkacki2016state}, as well as some moderate computational gains when compared to the standard PT.

  \begin{table}[tbh]
	\centering
				
\begin{tabular}{r|r|r}
	\cline{1-3}
	\multicolumn{1}{l}{} & \multicolumn{1}{l}{$\mathcal{C}^{1/2}_{k,\mathrm{Tempered}}$} & \multicolumn{1}{|l}{$\mathcal{C}^{1/2}_{k,\mathrm{RWM}}$} \\ \hline 
	$k=1$                & $\mathrm{Diag}(0.01,0.01,0.2,5)$                              & $\mathrm{Diag}(0.02,0.02,0.2,5)$                         \\ 
	$k=2$                & $\mathrm{Diag}(0.06,0.06,0.4,14)$                             & -                                                        \\ 
	$k=3$                & $\mathrm{Diag}(0.3,0.3,0.6,20)$                               & -                                                        \\ 
	$k=4$                & $\mathrm{Diag}(1,1,1,50)$                                     & -                                                        \\ 
\end{tabular}
	\caption{Step size of the RWM proposal distribution for the acoustic BIP experiment. Here $\mathrm{Diag}(d_1,d_2,\dots,d_N)$ is to be understood as the $N\times N$ diagonal matrix with entries $d_1,d_2,\dots,d_N$.   }

	\label{tab:size_ac}
\end{table}

\begin{table}[tbh]

	\centering
			\resizebox{\columnwidth}{!}{%
		\begin{tabular}{lcccccc}
			\hline
			& \multicolumn{2}{c}{Mean} & \multicolumn{2}{c}{Var} &  \multicolumn{2}{c}{Var$_\mathrm{RWM}/$Var}\\
			& $s_1$     & $s_2$          & $s_1$     & $s_2$   & $s_1$     & $s_2$ \\ \hline
			\multicolumn{1}{c}{} &             &            &               &              &                          \\
			RWM                  &1.33801    & 1.54293       & $9.86\times10^{-1}$    &   8.21$\times10^{-2}$  &1.000000&  1.000   \\
			PT              	    &  1.50121    & 1.00829        & $6.61\times 10^{-6}$    & $2.77\times 10 ^{-4}$    &149136.1& 296.2    \\
			PSDPT		         & 1.39775     & 1.23119       & $2.48\times10^{-1}$    & $6.54\times 10^{-2}$   &3.900000& 1.200    \\
			UGPT                & 1.50177     &1.00711        & $2.72\times10^{-6}$    & $2.38\times10^{-4}$    &361744.5& 345.0   \\
			WGPT               & 1.50174     & 1.00601        & $2.08\times10^{-6}$    & $1.46\times10^{-4}$    & 472133.2& 558.6     \\\hline
	\end{tabular}}

	\caption{Results for the acoustic BIP problem. Once again, we can see that both GPT algorithm provide an improvement over RWM, PT and PSDPT. The computational cost is comparable across all algorithms, given that the cost of each iteration is dominated by the cost of solving the underlying PDE. }
	\label{tab:res_ac}

\end{table}

\subsection{\rjp{High-dimensional acoustic wave inversion}}\label{ss:high_dim}

Lastly, we present a high-dimensional example for which we try to invert for the material properties $\beta^2$ in \eqref{Eq:acoustic}. 
For simplicity, we will consider fixed values of $\alpha=1,$  $s_1=1.5,$ and $s_2=1$.  In this case, we set $\beta^2=10+\hat{\beta}^2(x)$, where $\hat{\beta}(x)$ is taken  to be a realization of a random field discretized on a mesh of $N_x\times N_y$ triangular elements. This modeling choice ensures that $\beta^2$ is lower bounded. In this case, we will invert for the nodal values of (the finite element discretization of) $\hat{\beta}$, which will naturally result in a high-dimensional problem.  We remark that one is usually interested in including the randomness in $\beta^2$, instead of $\hat{\beta}$; however, we purposely choose to do so to induce an explicitly multi-modal posterior, and as such, to better illustrate the advantages of our proposed methods when sampling from these types of distributions. 

We begin by formalizing the finite-element discretization of the parameter space (see e.g., \cite{bui2016fem} for a more detailed discussion).

Let $\bar D=[0,3]\times[0,2],$ denote the physical space of the problem and let $V_h$ be a finite-dimensional subspace of $L_2(D)$ arising from a given finite element discretization. We write the finite element approximation $\hat{\beta}_h\in V_h$ of $\hat{\beta}$ as $$\hat{\beta}(x)\approx\hat{\beta}_h(x)=\sum_{n=1}^{N_\text{v}}b_n\phi_n(x),$$ where $\{\phi\}_{n=1}^{N_\text{v}}$ are the Lagrange basis functions corresponding to the nodal points $\{x_n\}_{n=1}^{N_\text{v}}$,  $(b_1,\dots,b_{N_\text{v}})^T=:\te\in\R^{N_\text{v}}$ is the set of nodal parameters and $N_\text{v}$ corresponds to the number of vertices used in the FE discretization. Thus, the problem of inferring the distribution of $\beta$ given some data $y$, can be understood as inferring the distribution of $\te$  given $y$. For our particular case, we will discretize $D$ using $28\times 28$ (non-overlapping) piece-wise linear finite elements, which results in $N_\text{v}=841$ and as such $\Theta=\R^{841}$. We consider a Gaussian prior $\mu_\mathrm{pr,\infty}=\mathcal{N}(0,\mathcal{A}^{-2})$, where $\mathcal{A}$ is  a differential operator acting on $L_2(D)$ of the form
 \begin{align}\label{eq:prec_op}
\mathcal{A}:=-a\nabla\cdot(H \nabla )+d I, \quad a,d>0,
\end{align}
together with Robin boundary conditions $\nabla (\cdot) \cdot \hat n +\sqrt{ad}(\cdot)=0$, where, following \cite{VillaPetraGhattas21},   $H$ is taken of the form  \begin{align}
	H:=\begin{pmatrix}
		e_1\sin^2(\ell)+e_2\cos^2(\ell)&&& (e_1-e_2)\sin(\ell)\cos(\ell)\\		
		(e_1-e_2)\sin(\ell)\cos(\ell)&&&e_1\cos^2(\ell)+e_2\sin^2(\ell)
	\end{pmatrix}.
\end{align}
Here $H$ models the spatial anisotropy of a Gaussian Random field sampled from $\mu_\mathrm{pr,\infty}$. 
It is known that for a two-dimensional (spatial) space, the covariance operator $\mathcal{A}^{-2}$ is symmetric and trace-class \cite{bui2016fem}, and as such, the  (infinite-dimensional) prior measure is well-defined. Thus, we set
\begin{align}
	\hat{\beta}(x)\sim{\mu}_\mathrm{pr,\infty},
\end{align}
which in turn induces the discretized prior:
 \begin{align}
 	\rjpr{\hat{\beta}_h}(x)\sim\mu_\mathrm{pr}:=\mathcal{N}(0,\mathcal{A}_h^{-2}),
\end{align}
where $\mathcal{A}_h^{-2}$ is a finite-element approximation of $\mathcal{A}$ using $28\times 28$ (non-overlapping) piece-wise linear finite elements. Samples from $\mu_\mathrm{pr}$ are obtained using the \texttt{FEniCS} package \cite{logg2012automated} and the \texttt{hIPPYlib} library \cite{VillaPetraGhattas21}.

\begin{figure}
	\centering
	\includegraphics[width=1.\linewidth]{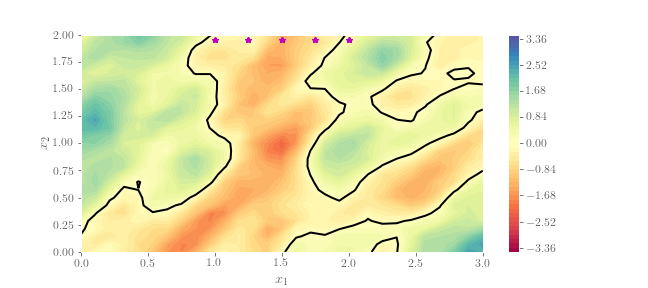}
	\caption{True field $\hat{\beta}^*_h(x)$. Notice the anisotropy on the field. The magenta points represent the receiver locations. The black line represents the zero-level set of the field. }
	\label{fig:truefield}
\end{figure}

We follow an approach similar to our previous example. We collect data $y\in Y$ by solving Equation \eqref{Eq:acoustic} with a force term given by \eqref{eq:force} and a true field $\hat{\beta}^*_h\sim\mu_\mathrm{pr}$ with $a=0.1$, $d=0.5$, $\ell=\pi/4$, $e_1=2$ and $e_2=0.5$. \rjpr{Such a realization of $\hat{\beta}^*_h$ is shown in Figure \ref{fig:truefield}.}

 Furthermore, data is observed at $N_R=5$ different receiver locations $R_1=(1.0,2.0),R_2=(1.25,2.0),R_3=(1.5,2.0),R_4=(1.75,2.0)$, and $R_5=(2.0,2.0)$  at $N_T=600$ equally-spaced time instants between $t=0$ and $t=0.6$. The data measured by each receiver is polluted by an (independent) additive Gaussian noise $\eta\sim\mathcal{N}(0,\sigma_\mathrm{noise}^2 I_{600\times 600})$, with $\sigma=0.021$, which corresponds to roughly a 0.5\% noise.  Thus, we have that $(Y,\lno \cdot\rno_Y)=(\R^{5\times 600},\lno \cdot\rno_\Sigma)$.  Similarly as in Section \ref{ss:exp_bip_ac}, the forward mapping operator $\mathcal{F}:\Theta\mapsto Y$ can be understood as the numerical solution of Equation \eqref{Eq:acoustic} evaluated at 600 discrete time instants at each of the 5 receiver locations. Numerical implementation follows a similar set-up as in Section \ref{ss:exp_bip_ac}, however,  for simplicity, we use $28\times 28$ triangular elements to approximate the forward operator $\eff$. The Bayesian inverse problem at hand can thus be understood as sampling from the posterior measure $\mu$, which has a Radon-Nikodym derivative with respect to the prior $\mu_\mathrm{pr}$ given by \begin{align}\label{Eq:ac_post2}
	\pi(\te)=\frac{\mathrm{d} \mu}{\mathrm{d} \mu_\mathrm{pr}}(\te)=\frac{1}{Z}\exp\left(-\frac{1}{2}\lno y-\mathcal{F}(\te)\rno^2_\Sigma\right).
\end{align}

\begin{figure}
	\centering
	\includegraphics[width=1\linewidth]{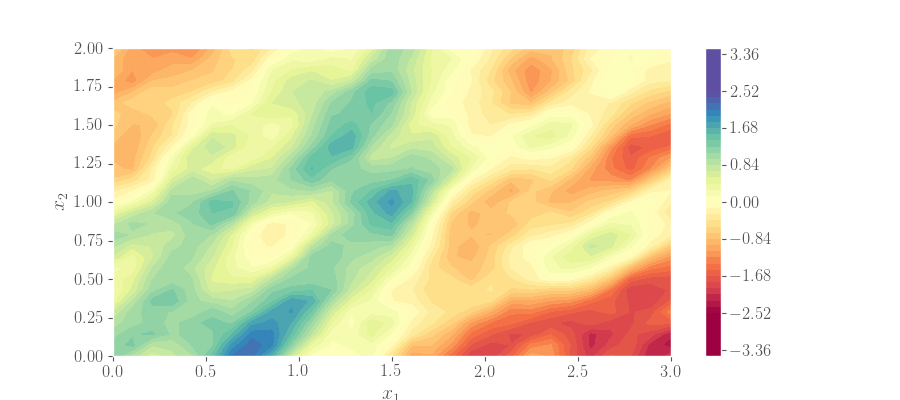}\\
	\includegraphics[width=1\linewidth]{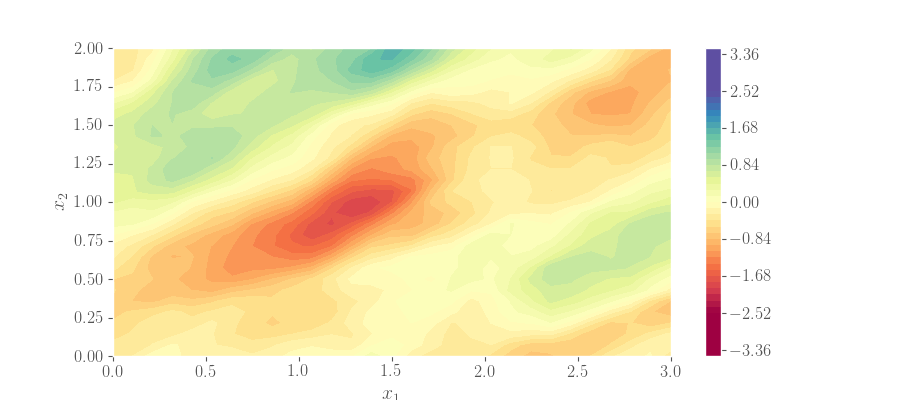}\\
	\includegraphics[width=1\linewidth]{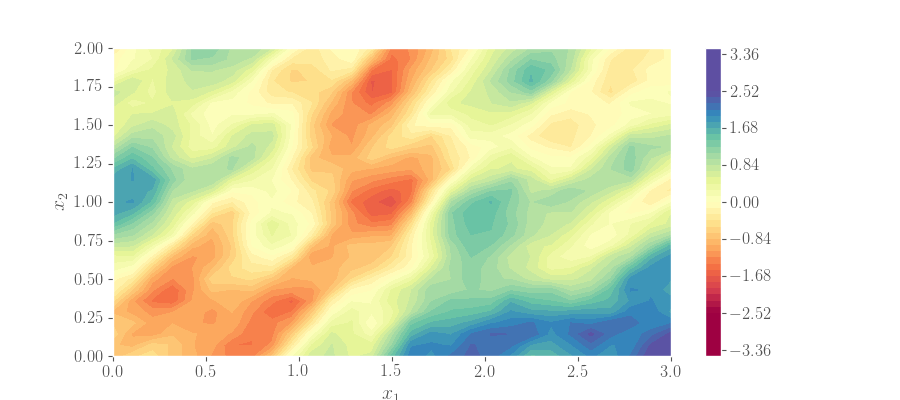}\\
	\includegraphics[width=1\linewidth]{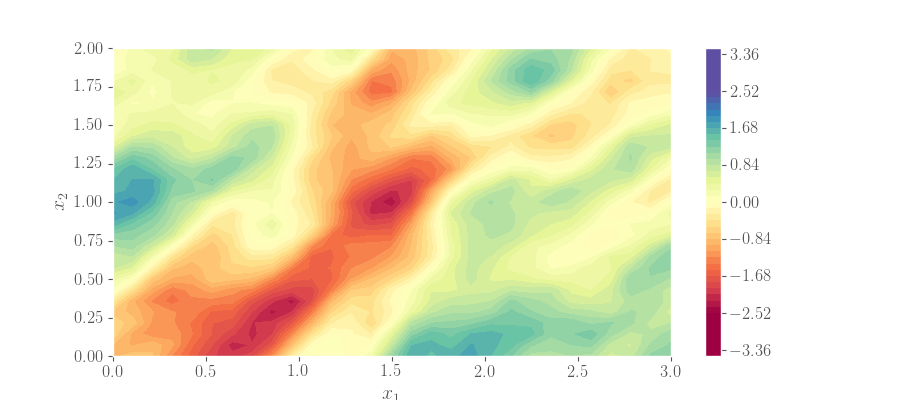}\\
	\caption{Posterior samples \rjpr{$\hat{\beta}_h\sim\mu$} obtained with the UW GPT algorithm. Notice the resemblance to Figure \ref{fig:truefield}.}
	\label{fig:post_samples}
\end{figure}

 The previous BIP has several difficulties; clearly, it is a high-dimensional posterior. Furthermore, just as in the previous example, the underlying mathematical model for the forward operator is a costly time-dependent PDE. By choosing to invert for \rjpr{ $\hat{\beta}_h\sim \mu_\mathrm{pr}$ (instead of $\beta_h^2\approx \beta^2$)}, and since $\mu_\mathrm{pr}$ is centered at zero, we induce a  multi-modal posterior, indeed, if the posterior concentrates around $\hat{\beta}^*_h$ it will also have peaks at any other $\hat{\beta}^j$ obtained by flipping the sign of $\hat{\beta}^*_h$ in a concentrated region separated by the zero level set of $\hat{\beta}^*_h$ (we identify 7 regions in Figure \ref{fig:truefield}). This can be seen in Figure \ref{fig:post_samples}, where we plot 4 \rjpr{posterior samples $\hat{\beta}_h\sim\mu$}. Notice the change in sign between some regions. Lastly, as a quantities of interest, we will consider $\mathcal{Q}_1=\int_D \exp(\hat{\beta}(x))\mathrm{d}x$ and $\mathcal{Q}_2=\exp(\hat{\beta}(1.5,1))$. We remark that, although these quantities of interest do not have any meaningful physical interpretation, they are, however, affected by the multi-modality of the posterior, and as such, well suited to exemplify the capabilities of our method.

Given the high-dimensionality of the posterior, we present a slightly different experimental setup in order to estimate $\E_{\mu}[\mathcal{Q}_i]\approx \hat{\mathcal{Q}}_i$, $i=1,2$. In particular, we will use the preconditioned Crank-Nicolson (pCN)  as a base method, instead of  RWM, for the transition kernel $p$. We compare the quality of our algorithms by examining the variance of the estimators $\hat {\mathcal{Q}}_i$ computed over $N_\text{runs}=50$ independent MCMC runs of each algorithm, with $K=4$ temperatures for the tempered algorithms  given by $T_1=1, T_2=4.57, T_3=20.89, T_4=100$. For the tempered algorithms, each estimator is obtained by running the inversion experiment for $N=4,800$ samples, discarding the first 20\%  of the samples (800) as a burn-in. For the untempered pCN algorithm, we run the inversion experiment for  $N_\mathrm{pCN}=K N=19, 200$ iterations, and discard the first 20\% of the  samples obtained (3840) as a burn-in.

Each individual chain is constructed using pCN proposals $q_{\mathrm{prop},k}(\te^{n}_k,\cdot)=\mathcal{N}(\sqrt{1-\rho_k^2}\te^{n}_k,\rho_k^2 \mathcal{A}_h^{-2})$, $k=1,2,3,4$, with $\rho_k$ described in Table \ref{tab:size_wave}. The simple, un-tempered pCN algorithm is run with a step size given by $\rho=\rho_1$. The values of $\rho_k$ are tuned in such a way that the acceptance rate of each chain is around 0.3 and are reported in Table \ref{tab:size_wave}.  The variance of the estimators obtained with each method is presented in Table \ref{tab:res_ac2}. Once again, even for this high-dimensional, highly multi-modal case, our proposed methods perform considerably better than the other algorithms.

\begin{table}[tbh]
	\centering
	\begin{tabular}{lrrrr}
		\hline
		&        $k=1$ &  $k=2$ &  $k=3$ &        $k=4$  \\
		\hline
		$\rho_{k}$     &  0.1 &  0.2 &  0.4 &  0.8 \\
		\hline
	\end{tabular}
	\caption{Values of $\rho_k$ for the pCN kernel for the high-dimensional wave inversion problem.   }
	
	\label{tab:size_wave}
\end{table}

\begin{table}[tbh]
	
	\centering
\small{
\begin{tabular}{lcccccc}
	\hline
	& \multicolumn{2}{c}{Mean} & \multicolumn{2}{c}{Var} &  \multicolumn{2}{c}{Var$_\mathrm{pCN}/$Var}\\
	& $\hat{\mathcal{Q}}_1$     & $\hat{\mathcal{Q}}_2$          & $\hat{\mathcal{Q}}_1$     & $\hat{\mathcal{Q}}_2$   & $\hat{\mathcal{Q}}_1$     & $\hat{\mathcal{Q}}_2$ \\ \hline
	\multicolumn{1}{c}{} &             &            &               &              &                          \\
	pCN                 &8.8665     & 1.5255       &  5.7362   &   0.6029  &1.00&  1.00   \\
	PT +pCN             	 &  8.7710    & 1.5311        & 1.3308    & 0.1380    &4.31& 4.36    \\
	PSDPT	+pCN	         & 8.5546     & 1.4453       &  2.1289    & 0.2666   &2.69& 2.26   \\
	UGPT    +pCN            & 8.7983     &1.4614        & 1.0543     &  0.1051    &5.49& 5.73   \\
	WGPT     +pCN          & 8.6464     & 1.4643       & 1.0126    	& 0.1016    & 5.74& 5.93     \\\hline
\end{tabular}	}

	\caption{Results for the high-dimensional acoustic BIP problem. As for the previous examples, The computational cost is comparable across all algorithms. }
	\label{tab:res_ac2}
	
\end{table}

\section{Conclusions and future work}\label{S:conclusions}
In the current work, we have proposed, implemented, and analyzed two  versions of the GPT, and applied these methods to a BIP context. We demonstrate that such algorithms produce reversible and geometrically-ergodic chains under relatively  mild conditions. As shown in Section \ref{S:Num_Exp}, such sampling algorithms provide an attractive alternative to the more standard Parallel Tempering when sampling from \textit{difficult} (i.e., multi-modal or concentrated around a manifold) posteriors. We remark that the framework considered here-in can be combined with other, more advanced MCMC algorithms, such as, e.g., the Metropolis-adjusted Langevin algorithm (MALA), or the Delayed Rejection Adaptive Metropolis (DRAM), for example \cite{haario2006dram}. 

We intend to carry out a number of future extensions of the work presented herein. One of our short-term goals is to extend the methodology developed in the current work to a Multi-level Markov Chain Monte Carlo context, as in \cite{dodwell2015hierarchical,madrigal2021analysis}  In addition, from a theoretical point of view, we would like to investigate the role that the number of chains and the choice of temperatures play on the convergence of the  GPT algorithm, as it has been done previously for Parallel Tempering in \cite{woodard2009conditions}. Improving on the estimates presented here would likely be the focus of future work. We also  believe that by excluding the identity permutation (i.e., $\mathrm{id}\notin S_K$) on the UGPT, one could obtain a swapping kernel which is \emph{better} in the so-called Peskun sense, see \cite{andrieu2009pseudo} for more details. We intend to carry further numerical experiments to better understand and compare swapping strategies.
Furthermore, from a computational perspective, given that the framework presented in this work is, in principle, dimension independent, the methods explored in this work can  also be combined with dimension-independent samplers such as the ones presented in \cite{beskos2017geometric,cui2016dimension}, thus providing a sampling algorithm robust to both multi-modality and large dimensionality of the parameter space. Given the additional computational cost of these methods, a non-trivial coupling of GPT and these methods needs to be devised. 
Lastly, we aim at applying the methods developed in the current work to more computationally challenging BIP, in particular those arising in seismology and seismic source inversion, where it is not uncommon to find multi-modal posterior distributions when inverting for a point source.  
\begin{appendix}
	\section*{Appendix}\label{appn} 

	\section{Proof of Proposition \ref{proposition:bound_var}}\label{AP:proof_bound_cov}
\begin{proof}This proof is partially based on the proof of Theorem 1.2 in \cite{madras2002markov}.
	Let $\teb,\yeb\in \Theta^K$ and define $\fbs:=\mub_\sigma(f)$. The right-most inequality follows from the fact that 
	\begin{align}
	1&=\V_{\mubis}[f]=\int_{\Theta^K} f(\teb)^2\mubis(\mathrm{d}\teb)\\&=\frac{1}{|S_K|}\sum_{\sigma\in S_K}\int_{\Theta^K}f^2(\teb)\mub_\sigma(\mathrm{d}\teb)\\&=\frac{1}{|S_K|}\sum_{\sigma\in S_K} \left(\V_{\mub_\sigma}[f]+\bar{f}_\sigma^2\right)\geq\frac{1}{|S_K|}\sum_{\sigma\in S_K} \V_{\mub_\sigma}[f]
	\end{align}
We follow a procedure similar to the proof of \cite[Theorem 1.2]{madras2002markov} for the lower bound on the variance. We introduce an ordering on $S_K={\sigma_1,\sigma_2,\dots,\sigma_{|S_K|}},$  define the matrix $C\in \R^{{|S_K|}\times {|S_K|}}$ as the matrix with entries \begin{align}\label{Eq:cij}
C_{ij}=\int_{\Theta^K}\int_{\Theta^K} (f(\teb)-f(\yeb))^2\mub_{\sigma_i}(\mathrm{d}\teb)\mub_{\sigma_j}(\mathrm{d}\yeb), 
\end{align}
where $C_{jj}=2\V_{\mub_{\sigma_j}}[f]$ and
\begin{align}
2=2\V_{\mubis}[f]&=\int_{\Theta^K}\int_{\Theta^K}(f(\teb)-f(\yeb))^2\left(\frac{1}{{|S_K|}}\sum_{i=1}^{|S_K|}\mub_{\sigma_i}(\mathrm{d}\teb)\right)\\&\times\left(\frac{1}{{|S_K|}}\sum_{j=1}^{|S_K|}\mub_{\sigma_j}(\mathrm{d}\yeb)\right)\\
&=\sum_{i,j}\frac{1}{{|S_K|^2}}C_{ij} \label{Eq:cov_mat}.
\end{align}
We thus aim at finding an upper bound of Equation (\ref{Eq:cov_mat}) in terms of $({|S_K|})^{-1}\sum_{\sigma\in S_K}\V_\sigma[f]$. 

By assumption \ref{ass:3}, for any $\sigma_i,\sigma_j\in S_K$ the densities $\pib_{\sigma_i},\pib_{\sigma_j}$ of $\mub_{\sigma_i},\mub_{\sigma_j}$ (with respect to $\mub_\mathrm{pr}$) have an overlap $\Lambda_{\sigma_i,\sigma_j}>0$. For brevity, in the following we use the shorthand notation $\Lambda_{i,j}$ for $\Lambda_{\sigma_i,\sigma_j}$
Thus, we can find densities $$\etab_{ij}:=\dij^{-1}\underset{\teb\in \Theta^K}{\min}\{ \pib_{\sigma_i}(\teb),\pib_{\sigma_j}(\teb) \},\fib,\psib$$  such that 
$\pib_{\sigma_i}=\dij\etab_{ij}+(1-\dij)\fib,$ and $\pib_{\sigma_j}=\dij\etabij +(1-\dij)\psib$. Thus, integrating over $\Theta^K,$ we get for the diagonal entries of the $C$ matrix:

\begin{align}
&C_{ii}=2\V\mub_{\sigma_i}[f]\\&=\iint(f(\teb)-f(\yeb))^2\left(\dij\etab_{ij}(\teb)+(1-\dij)\fib(\teb)\right)\\&\times\left(\dij\etabij(\yeb) +(1-\dij)\fib(\yeb)\right)\mub_\mathrm{pr}(\mathrm{d}\teb)\mub_\mathrm{pr}(\mathrm{d}\yeb)\\
&=\iint(f(\teb)-f(\yeb))^2\dij^2\etabij(\teb)\etabij(\yeb)\mub_\mathrm{pr}(\mathrm{d}\teb)\mub_\mathrm{pr}(\mathrm{d}\yeb)\\
&+\iint(f(\teb)-f(\yeb))^2\dij(1-\dij)\fib(\yeb)\etabij(\teb)\mub_\mathrm{pr}(\mathrm{d}\teb)\mub_\mathrm{pr}(\mathrm{d}\yeb)\\
&+\iint(f(\teb)-f(\yeb))^2\dij(1-\dij)\fib(\teb)\etabij(\yeb)\mub_\mathrm{pr}(\mathrm{d}\teb)\mub_\mathrm{pr}(\mathrm{d}\yeb)\\&+\iint(f(\teb)-f(\yeb))^2(1-\dij)^2\fib(\yeb)\fib(\teb)\mub_\mathrm{pr}(\mathrm{d}\teb)\mub_\mathrm{pr}(\mathrm{d}\yeb)\\
&=2\dij^2\V_{\etabij}[f]+2(1-\dij)^2\V_{\fib}[f]+2\dij(1-\dij)\\
&\times \iint(f(\teb)-f(\yeb))^2\etabij(\teb)\fib(\teb)\mub_\mathrm{pr}(\mathrm{d}\teb)\mub_\mathrm{pr}(\mathrm{d}\yeb) \label{Eq:2vsigma}.
\end{align}
Notice that equation (\ref{Eq:2vsigma}) implies that 
\begin{align}\label{Eq:bound_vj_veta}
&\iint(f(\teb)-f(\yeb))^2\etabij(\teb)\fib(\teb)\mub_\mathrm{pr}(\mathrm{d}\teb)\mub_\mathrm{pr}(\mathrm{d}\yeb)\\
&\leq \frac{\V_{\mub_{\sigma_i}}[f]-\dij^2\V_{\etabij}[f]}{\dij(1-\dij)}.
\end{align}
As for the non-diagonal entries of $C$, we have
\begin{align}\label{Eq:cij_1}
&C_{ij}=\iint(f(\teb)-f(\yeb))^2\left[\dij\etab_{ij}(\teb)\right. \\&+ \left.(1-\dij)\fib(\teb)\right](\dij\etabij(\yeb)\\ &+(1-\dij)\psib(\yeb))\mub_\mathrm{pr}(\mathrm{d}\teb)\mub_\mathrm{pr}(\mathrm{d}\yeb)\\
&=2\dij^2\V_{\etabij}[f]\\&+(1-\dij)^2\iint(f(\teb)-f(\yeb))^2\fib(\teb)\psib(\yeb)\mub_\mathrm{pr}(\mathrm{d}\teb)\mub_\mathrm{pr}(\mathrm{d}\yeb)\\&+\dij(1-\dij)\iint(f(\teb)-f(\yeb))^2\\ &\times\left(\etabij(\teb)\psib(\yeb)+\etabij(\yeb)\fib(\teb)\right)\mub_\mathrm{pr}(\mathrm{d}\teb)\mub_\mathrm{pr}(\mathrm{d}\yeb). 
\end{align}
We can bound the second term in the previous expression using Cauchy-Schwarz. Let $\bm z \in \Theta^K$. Then, 
\begin{align}
&\iint(f(\teb)-f(\yeb))^2\fib(\teb)\psib(\yeb)\mub_\mathrm{pr}(\mathrm{d}\teb)\mub_\mathrm{pr}(\mathrm{d}\yeb)\\
&=\iiint(f(\teb)-f(\bm z)+f(\bm z)-f(\yeb))^2\fib(\teb)\psib(\yeb)\etabij(\bm z)\\&\times\mub_\mathrm{pr}(\mathrm{d}\teb)\mub_\mathrm{pr}(\mathrm{d}\yeb)\mub_\mathrm{pr}(\mathrm{d}\bm z)\\
&\leq 2\iiint\left( \left(f(\teb)-f(\bm z)\right)^2+(f(\bm z)-f(\yeb))^2\right)\fib(\teb)\psib(\yeb)\etabij(\bm z)\\&\times\mub_\mathrm{pr}(\mathrm{d}\teb)\mub_\mathrm{pr}(\mathrm{d}\yeb)\mub_\mathrm{pr}(\mathrm{d}\bm z)\\
&=2\iint(f(\teb)-f(\bm z))^2\fib(\teb)\etabij(\bm z)\mub_\mathrm{pr}(\mathrm{d}\teb)\mub_\mathrm{pr}(\mathrm{d}\bm z)\\
&+2\iint(f(\yeb)-f(\bm z))^2\psib(\yeb)\etabij(\bm z)\mub_\mathrm{pr}(\mathrm{d}\yeb)\mub_\mathrm{pr}(\mathrm{d}\bm z). \label{Eq:bound_with_cs}
\end{align}
Thus, from equations (\ref{Eq:bound_vj_veta}), (\ref{Eq:cij_1}), and (\ref{Eq:bound_with_cs}) we get
\begin{align}
C_{ij}&\leq 2\dij^2\V_{\etabij}[f]+ (2(1-\dij)^2+\dij(1-\dij))\\&\times\left(\iint(f(\teb)-f(\bm y))^2\left(\etabij(\teb)\bm \psi_j(\bm y) \right. \right. \\&\left.\left.+\etabij(\bm y)\psi_i(\bm \teb)\right)\mub_\mathrm{pr}(\mathrm{d}\teb)\mub_\mathrm{pr}(\mathrm{d}\bm y) \right)\\
&=2\dij^2 \V_{\etabij}[f]+(2-\dij)(1-\dij)\\&\times\frac{\left(\V_{\mub_{\sigma_i}}[f]-\dij^2\V_{\etabij}[f]+\V_{\mub_{\sigma_j}}[f]-\dij^2\V_{\etabij}[f]\right)}{\dij(1-\dij)}\\
&= \frac{2-\dij}{\dij}\left(V_{\mub_{\sigma_i}}[f]+V_{\mub_{\sigma_j}}[f]\right)-4\dij(1-\dij)\V_{\etabij}[f]\\
&\leq\frac{2-\dij}{\dij}\left(V_{\mub_{\sigma_i}}[f]+V_{\mub_{\sigma_j}}[f]\right)\label{Eq:bound_for_cs},
\end{align}
since $\dij\in(0,1)\ \forall i,j$. Finally, from equations (\ref{Eq:cov_mat}) and (\ref{Eq:bound_for_cs}) we get that \begin{align}
1&=V_{\mubis}[f]=\frac{1}{2}\sum_{i,j}\frac{1}{{|S_K|^2}}C_{ij}\\&\leq\frac{1}{2}\frac{1}{{|S_K|^2}}\sum_{i,j=1}^{|S_K|}\frac{2-\dij}{\dij}\left(V_{\mub_{\sigma_j}}[f]+V_{\mub_{\sigma_j}}[f]\right)\\ &\leq\frac{2-\Lambda_m}{\Lambda_m}\left(\frac{1}{{|S_K|}}\sum_{i=1}^{|S_K|}\V_{\mub_{\sigma_i}}[f]\right),
\end{align}
with $\Lambda_m:=\underset{i,j=1,2,\dots,|S_K|}{\min\{ \dij\}}>0$, and $\Lambda_{i,j}$ as in Assumption \ref{ass:3}. Notice that we have used (\ref{Eq:bound_for_cs}) for the first inequality, including the case $i=j$, in the previous equation. This in turn yields the lower bound  
\begin{align}
0<\frac{\Lambda_m}{2-\Lambda_m}\leq\left(\frac{1}{{|S_K|}}\sum_{i\in S_K}\V_{\mub_{_i}}[f]\right).
\end{align}
\end{proof}

\end{appendix}
%
%
%
%
\section*{Acknowledgements}
\rjp{We would like to thank the anonymous reviewers for helpful suggestions that significantly improved this work.} This publication was supported by funding from King Abdullah University of Science and Technology (KAUST) Office of Sponsored Research (OSR) under award numbers URF/1/2281-01-01 and URF/1/2584-01-01 in the KAUST Competitive Research Grants Programs- Round 3 and 4, respectively, and the Alexander von Humboldt Foundation.  Jonas Latz acknowledges support by the Deutsche Forschungsgemeinschaft (DFG) through the TUM International
Graduate School of Science and Engineering (IGSSE) within the project 10.02 BAYES. Juan P. Madrigal-Cianci and Fabio Nobile also acknowledge support from  the {Swiss Data Science Center (SDSC) Grant p18-09.


%
%



\bibliographystyle{spmpsci} 
\bibliography{Generalized_PT_STATCO}       

\end{document}